\documentclass[letterpaper]{article}

\author{Benjamin Antieau}
\usepackage[letterpaper,left=1.5in,right=1.5in,top=1.5in,bottom=1.5in]{geometry}

\newcommand{\myauthor}{Benjamin Antieau}
\newcommand{\mytitle}{\'Etale twists in noncommutative algebraic geometry and the twisted Brauer space}
\newcommand{\pdftitle}{\mytitle}

\title{\'Etale twists in noncommutative algebraic geometry and the twisted Brauer space}

\usepackage[pdfstartview=FitH,
            pdfauthor={\myauthor},
            pdftitle={\pdftitle},
            colorlinks,
            linkcolor=reference,
            citecolor=citation,
            urlcolor=e-mail,
            ]{hyperref}

\usepackage{mathrsfs}
\usepackage[mathscr]{euscript}
\usepackage{amsmath}
\usepackage{amscd}
\usepackage{amsbsy}
\usepackage{amssymb}
\usepackage{microtype}
\usepackage{times}
\usepackage{enumerate}
\usepackage[all,cmtip]{xy}
\usepackage{tikz}
\usetikzlibrary{matrix,arrows}
\usepackage{dsfont}
\usepackage[abbrev,lite]{amsrefs}
\usepackage{amsthm}

\usepackage{color}
\definecolor{todo}{rgb}{1,0,0}
\definecolor{conditional}{rgb}{0,1,0}
\definecolor{e-mail}{rgb}{0,.40,.80}
\definecolor{reference}{rgb}{.20,.60,.22}
\definecolor{mrnumber}{rgb}{.80,.40,0}
\definecolor{citation}{rgb}{0,.40,.80}

\DeclareMathOperator{\Ho}{Ho}

\newcommand{\perf}{\mathrm{perf}}

\DeclareMathOperator{\id}{id}

\DeclareMathOperator{\op}{op}

\DeclareMathOperator{\im}{im}

\newcommand{\desc}{\mathrm{desc}}

\newcommand{\rwe}{\tilde{\rightarrow}}

\newcommand{\we}{\simeq}
\newcommand{\iso}{\cong}

\newcommand{\Gm}{\mathds{G}_{m}}

\newcommand{\Gal}{\mathrm{Gal}}

\newcommand{\ind}{\mathrm{ind}}

\newcommand{\tors}{\mathrm{tors}}

\newcommand{\qc}{\mathrm{qc}}
\newcommand{\cont}{\mathrm{cont}}
\newcommand{\ccc}{\mathrm{c}}
\newcommand{\sep}{\mathrm{sep}}

\DeclareMathOperator{\THH}{HH}
\newcommand{\THHB}{\mathbf{HH}}
\DeclareMathOperator{\K}{\mathrm{K}}

\newcommand{\Mod}{\mathrm{Mod}}
\newcommand{\Perf}{\mathrm{Perf}}

\newcommand{\Shv}{\mathrm{Shv}}

\newcommand{\Cat}{\mathrm{Cat}}

\DeclareMathOperator{\Spec}{Spec}

\DeclareMathOperator{\PGL}{PGL}
\newcommand{\ShPGL}{\mathbf{PGL}}
\newcommand{\ShPSO}{\mathbf{PSO}}

\DeclareMathOperator{\SL}{SL}

\DeclareMathOperator{\Pic}{Pic}

\DeclareMathOperator{\Hoh}{H}
\DeclareMathOperator{\Eoh}{E}

\DeclareMathOperator{\Roh}{R}

\newcommand{\Hom}{\mathrm{Hom}}
\newcommand{\Map}{\mathrm{Map}} 
\newcommand{\Fun}{\mathrm{Fun}}
\newcommand{\map}{\mathrm{map}} 

\DeclareMathOperator{\End}{End}
\DeclareMathOperator{\Aut}{Aut}

\newcommand{\ShPic}{\mathbf{Pic}}

\newcommand{\ShPr}{\mathbf{Pr}}

\newcommand{\ShBr}{\mathbf{Br}}

\newcommand{\ShMor}{\mathbf{Mor}}

\newcommand{\Shaut}{\mathbf{aut}}
\newcommand{\ShAut}{\mathbf{Aut}}
\newcommand{\ShB}{\mathbf{B}} 

\newcommand{\StCat}{\mathscr{C}\mathrm{at}}

\newcommand{\StMod}{\mathscr{M}\mathrm{od}}

\newcommand{\et}{\mathrm{\acute{e}t}}

\DeclareMathOperator{\Br}{Br}

\newcommand{\Lrm}{\mathrm{L}}

\newcommand{\Mrm}{\mathrm{M}}

\newcommand{\Hrm}{\mathrm{H}\,}
\newcommand{\Nrm}{\mathrm{N}}

\newcommand{\Drm}{\mathrm{D}}

\newcommand{\Fscr}{\mathscr{F}}
\newcommand{\Oscr}{\mathscr{O}}

\newcommand{\Ascr}{\mathscr{A}}
\newcommand{\Bscr}{\mathscr{B}}

\newcommand{\HH}{\mathds{H}}
\newcommand{\CC}{\mathds{C}}
\newcommand{\RR}{\mathds{R}}
\newcommand{\QQ}{\mathds{Q}}
\newcommand{\ZZ}{\mathds{Z}}
\newcommand{\EE}{\mathds{E}}

\newcommand{\PP}{\mathds{P}}

\newcommand{\NP}{\mathds{NP}}

\theoremstyle{plain}
\newtheorem{theorem}{Theorem}[section]
\newtheorem{lemma}[theorem]{Lemma}
\newtheorem{proposition}[theorem]{Proposition}

\newtheorem{conjecture}[theorem]{Conjecture}
\newtheorem{corollary}[theorem]{Corollary}
\newtheorem{scholium}[theorem]{Scholium}

\newtheorem{problem}[theorem]{Problem}

\theoremstyle{definition}
\newtheorem{definition}[theorem]{Definition}

\newtheorem{example}[theorem]{Example}
\newtheorem{question}[theorem]{Question}

\newtheorem{remark}[theorem]{Remark}

\let\oldmarginpar\marginpar
\renewcommand\marginpar[1]{\-\oldmarginpar[\raggedleft\footnotesize #1]%
{\raggedright\footnotesize #1}}

\begin{document}
\maketitle

\begin{abstract}

    \noindent
    This paper studies \'etale twists of derived categories of schemes and associative
    algebras. A general method, based on a new construction called the twisted Brauer space,
    is given for classifying \'etale twists, and a complete classification is carried out
    for genus $0$ curves, quadrics, and noncommutative projective spaces. A partial
    classification is given for curves of higher genus.
    The techniques build upon my recent
    work with David Gepner on the Brauer groups of commutative ring spectra.


    \paragraph{Mathematics Subject Classification 2010}
    Primary:
    \href{http://www.ams.org/mathscinet/msc/msc2010.html?t=14Fxx&btn=Current}{14F22},
    \href{http://www.ams.org/mathscinet/msc/msc2010.html?t=18Exx&btn=Current}{18E30}.
    Secondary:
    \href{http://www.ams.org/mathscinet/msc/msc2010.html?t=14Dxx&btn=Current}{14D20},
    \href{http://www.ams.org/mathscinet/msc/msc2010.html?t=14Fxx&btn=Current}{14F05},
    \href{http://www.ams.org/mathscinet/msc/msc2010.html?t=16Exx&btn=Current}{16E40}.

\end{abstract}

\tableofcontents

%
%
%
%

\section{Introduction}

\subsection{An example}\label{sec:example}

The purpose of this paper is to create a formalism for answering questions of the following
kind. Suppose that $X$ is a variety over a field $k$. How can one classify the $k$-linear derived
categories $\Drm$ such that $\Drm_{\overline{k}}\we\Drm^b(X_{\overline{k}})$? For the
purposes of the following example, please take this problem at face value and believe that
there is a good notion of such ``derived categories'' $\Drm$ together with a way to tensor
with $\overline{k}$. This will all be explained later in the introduction and in the rest of
the paper.

Allow me to begin the paper with a motivating example.
Let $\Br^{\PP^1}(\RR)$ denote the set of (derived equivalence classes of) $\RR$-linear derived categories $\Drm$
such that $\Drm_\CC\we\Drm^b(\PP^1_{\CC})$. Thus, $\Br^{\PP^1}(\RR)$ classifies derived categories that are \'etale locally
equivalent to the derived category of $\PP^1$. Its objects can be viewed as noncommutative
\'etale twists of the projective line. I call $\Br^{\PP^1}(\RR)$ the $\PP^1$-twisted Brauer
set of $\RR$. It is a pointed set, where the point is the category $\Drm^b(\PP^1_{\RR})$.

Consider the real path algebra $\RR Q$, where $Q$ is the quiver
$\bullet\rightrightarrows\bullet$. It is a result of Be{\u\i}linson~\cite{beilinson} that
$\Drm^b(\PP^1_{\RR})$ and $\Drm^b(\RR Q)$ are equivalent as $\RR$-linear triangulated categories.

There are two obvious ways to construct elements in $\Br^{\PP^1}(\RR)$. First, one can
tensor any given element with the quaternion algebra $\HH$. For instance, $\HH\otimes_\RR\RR
Q$ gives another $\RR$-algebra which becomes Morita equivalent to $\CC Q$ over $\CC$. Indeed,
$\CC\otimes_\RR\RR Q\iso\CC Q$, while $\CC\otimes_\RR\left(\HH\otimes_\RR\RR
Q\right)\iso\Mrm_2(\CC)\otimes_\CC\CC Q\iso\Mrm_2(\CC Q)$. Thus, $\Drm^b(\HH\otimes_\RR\RR
Q)$ is an element of $\Br^{\PP^1}(\RR)$. This derived category has a more geometric
interpretation: it is the derived category $\Drm^b(\PP^1_\RR,\alpha)$ of $\alpha$-twisted coherent sheaves on $\PP^1$,
where $\alpha$ is the class in $\Br(\PP^1_\RR)$
pulled back from $\HH$. Equivalently, $\Drm^b(\PP^1_\RR,\alpha)$ is the derived category of quaternionic vector
bundles on $\PP^1_\RR$. Since $\Br(\RR)=\ZZ/2\cdot\HH$,
there are no further iterations of the construction.

The algebras $\RR Q$ and
$\HH\otimes_\RR\RR Q$ represent distinct elements in the pointed set $\Br^{\PP^1}(\RR)$.
The easiest way to see this is via algebraic $K$-theory. The $K$-theory of $\PP^1$ (or
equivalently of $\RR Q$)
is $\K_*(\RR)\oplus\K_*(\RR)$ by Quillen's computation~\cite{quillen}*{Theorem 8.2.1}, while
the $K$-theory of $\HH\otimes_\RR\RR Q$ is
$\K_*(\HH)\oplus\K_*(\HH)$. The torsion part of $\K_1(\RR)\iso\RR^\times$ is $\ZZ/2$, while the
torsion part of $\K_1(\HH)\iso\HH^\times/[\HH^\times,\HH^\times]$ is $0$, where
$[\HH^\times,\HH^\times]$ is the commutator subgroup of $\HH^\times$. The point is that the reduced norm
$\K_1(\HH)\rightarrow\K_1(\RR)$ is injective by the theorem of Wang~\cite{wang}. But, clearly,
$-1\in\RR^*$ cannot be the reduced norm of a quaternion. Thus, $\HH\otimes_\RR\RR Q$ and
$\RR Q$ are not derived Morita equivalent.

The second obvious way to construct elements in $\Br^{\PP^1}(\RR)$ is to
look at another variety over $\Spec\RR$ that becomes isomorphic to $\PP^1$
over $\Spec\CC$. Up to isomorphism, there is only one such variety, which is the genus $0$
curve $C$ cut out by $x^2+y^2+z^2=0$ in $\PP^2_\RR$. Since this curve does not have an $\RR$-point,
it is not the projective line, but it becomes isomorphic to $\PP^1$ over $\CC$. Thus
$\Drm^b(C)$ represents another point of $\Br^{\PP^1}(\RR)$.
Interestingly, in this case, considering $\alpha$-twisted sheaves gives nothing new. Because
$C$ is the Severi-Brauer variety of $\HH$, the pullback of $\HH$ to $C$ has zero Brauer
class. Thus, $\Drm^b(C,\alpha)\we\Drm^b(C)$. To see that $\Drm^b(C)$ is distinct from
either of the module categories from the previous paragraph, note that its
$K$-theory is isomorphic to $\K_*(\RR)\oplus\K_*(\HH)$ by Quillen's computation of the
$K$-theory of Severi-Brauer varieties~\cite{quillen}*{Theorem 8.4.1},
and this is different from either of the other $K$-theories, by consideration of torsion in
degree $1$.

Thus, there are at least $3$ elements of $\Br^{\PP^1}(\RR)$, and there is an action
on these elements by $\Br(\RR)$, which is described above. The main point of this paper is to
develop methods that will allow a precise formulation of the problems of the type posed in the example,
and to give a computational tool for solving these problems, which I apply in many cases.
In particular, in Section~\ref{sec:curvesoverr} this
computational tool will be used to show that there are no other elements in
$\Br^{\PP^1}(\RR)$ besides those described already.

Every element of $\Br^{\PP^1}(\RR)$ is
represented by a category of modules over an associative algebra. This has already been remarked upon for
$\Drm^b(\PP^1_\RR)$ and $\Drm^b(\PP^1_\RR,\alpha)$. For the genus $0$ curve $C$, there is an
equivalence $\Drm^b(C)\we\Drm^b(A)$, where $A$ is the path algebra of the modulated quiver
$\RR:\bullet\rightrightarrows\bullet:\HH$. Modulated quivers were used to classify finite
dimensional hereditary algebras of finite representation type. For details, see
Dlab-Ringel~\cite{dlab-ringel}, where they are called species.

\subsection{Overview}\label{sec:overview}

The noncommutative algebraic geometry of the title is what Ginzburg has called noncommutative algebraic geometry ``in the
large,'' where one replaces schemes with derived categories of sheaves and isomorphisms with
derived Morita equivalences. This form of noncommutative algebraic geometry, which began
with the work of Be{\u\i}linson~\cite{beilinson}, has been reinforced by ideas
originating in string theory, where two varieties with equivalent derived
categories should describe the same physical theory. The mathematical theory has been pursued by
Bondal, Ginzburg, Kontsevich, Orlov, Rosenberg, and van den
Bergh to name just a few.
See~\cite{bondal-orlov,bondal-orlov-coherent,bondal-vandenbergh,ginzburg,kontsevich-rosenberg,vandenberg-blowing}.
Thus, if $A$ is an associative algebra, the derived category of $A$-modules $\Drm(A)$ is
viewed as a geometric object. Noncommutative algebraic geometry in the large is distinct
from both noncommutative algebraic geometry in the small and derived algebraic geometry. The
former is about noncommutative deformations of commutative rings and is modeled on the
coordinate algebras that arise in quantum mechanics. Derived algebraic geometry on the other
hand replaces ordinary commutative rings with ``derived'' commutative rings, which are
either simplicial commutative rings, commutative dg algebras, or commutative ring spectra.

This viewpoint is motivic in the sense that many classical
motivic invariants, such as Hochschild homology and $K$-theory, depend only on the derived
category.

To formulate this kind of geometry correctly requires a more flexible framework than
simply triangulated categories. Thus, $\Drm(A)$ is replaced by $\Mod_A$, the stable
$\infty$-category of right $A$-modules, and the triangulated category of complexes of
$\Oscr_X$-modules is replaced by $\Mod_X$, a stable $\infty$-categorical model for
$\Drm_{\qc}(X)$. Another option would be to use dg enhancements or $A_\infty$-categories.
Stable $\infty$-categories include all of these examples, and have the added benefit that,
for instance, one can do geometry over the sphere spectrum.

Since I am interested in developing a theory that works over the sphere, my commutative
rings will be connective commutative ($\EE_{\infty}$-)ring spectra, and my associative rings
will be $A_{\infty}$-ring spectra. The reader will lose little
in thinking of ordinary commutative rings and associative dg algebras. But, in any case, a
module over a ring $A$, even an ordinary associative ring, means an $A_\infty$-module. So,
over an ordinary associative ring, modules are really complexes of ordinary $A$-modules.

Recall that a compact object in a stable $\infty$-category $\Mrm$ is an object $x$ such that
the mapping space functor $\map_\Mrm(x,-)$ commutes with filtered colimits. This is the
appropriate generalization of compactness in triangulated categories having all coproducts, where a compact object
$x$ is one where taking maps out commutes with coproducts.
Compact objects are the
cornerstone of noncommutative algebraic geometry. When $X$ is a quasi-compact and
quasi-separated scheme, the compact objects of $\Mod_X$ are precisely the perfect complexes,
which are complexes of $\Oscr_X$-modules locally quasi-isomorphic to bounded complexes of
finitely generated vector bundles.

The compact objects are fundamental to derived Morita theory. If $A$ and $B$ are two
associative algebras, then to give an equivalence $F:\Mod_A\rwe\Mod_B$ is to give a
compact right $B$-module $F(A)$ such that $F(A)$ generates $\Mod_B$ and $\End_B(F(A))\we A$. Note
that I will use derived equivalence for any equivalence between stable $\infty$-categories,
and not for a triangulated equivalence $\Drm(A)\we\Drm(B)$, although once a functor
$F:\Mod_A\rightarrow\Mod_B$ is given, the property of it being an equivalence can be
detected on the homotopy categories.

When $X$ is a reasonable scheme (quasi-compact and quasi-separated), Bondal and
van den Bergh~\cite{bondal-vandenbergh}
showed that there is a single perfect complex $E$ that generates the entire derived category
$\Drm_{\qc}(X)$. Thus, derived Morita theory says that at the level of $\infty$-categories, there is
an equivalence $\Mod_X\we\Mod_A$, where $A=\End_{\Mod_X}(E)^{\op}$ is the derived endomorphism
algebra spectrum of $E$. The example of Beilinson's, that $\Drm^b(\PP^1)\we\Drm^b(\RR Q)$,
from the previous section is an especially nice example of this phenomenon. In
particular, the algebra $A$ is typically truly an $A_{\infty}$-algebra, and is not derived Morita equivalent to
any ordinary associative algebra. Bondal and ven den Bergh's theorem justifies the term
noncommutative algebraic geometry. Almost \emph{every} derived category that arises in ordinary
algebraic geometry is the module category for an $A_{\infty}$-algebra, or is built from such
a category.

Therefore, from the perspective of noncommutative algebraic geometry, derived categories of
algebras are a natural generalization of derived categories of schemes. Thus, the
first question to ask is when two algebras or schemes give rise to the same noncommutative
geometric
object. For algebras, the answer, abstractly, is the subject of derived Morita theory,
which goes back to Cline-Parshall-Scott~\cite{cline-parshall-scott}, Happel~\cite{happel}, and Rickard~\cite{rickard}, and has been developed by many
people for use in the study of finite-dimensional associative algebras and in block theory
for modular representation theory. In the dg setting, Keller~\cite{keller} and
To\"en~\cite{toen-morita} have worked out the
theory very nicely. For ring spectra, the theory follows from work of Schwede and
Shipley~\cite{schwede-shipley}. The problem of when two varieties $X$ and $Y$ are derived equivalent has
been the subject of a great deal of research by Bondal, Bridgeland, Huybrechts, Kawamata,
Orlov, Stellari, van den Bergh, and many, many others. For a comprehensive introduction to
the subject and the literature, see~\cite{huybrechts}.

Now that there is an excellent categorical framework for studying derived equivalences, and
since the work of many authors has provided a clear picture of when to expect derived
equivalences, the follow-up question I want to ask in this paper is: is it possible to classify when
two algebras $A$ and $B$, say over a field $k$, represent the same geometric object over $\overline{k}$?
In fact, in general, it is better to ask for a finite separable extension $l/k$ such that
$A_l$ and $B_l$ are derived Morita equivalent. The analogous question for potentially
infinite or
inseparable extensions is considered in a special case in Section~\ref{sec:moritadescent}.

\begin{problem}\label{prob:classification}
    Let $A$ be an $A_{\infty}$-algebra over $k$. Classify, up to derived equivalence over
    $k$, all $A_{\infty}$-algebras $B$ such that
    $A\otimes_k k^{\sep}$ and $B\otimes_k k^{\sep}$ are derived Morita equivalent.
\end{problem}

When $\Mod_A\we\Mod_X$, the algebras $B$ should be viewed as noncommutative \'etale twists of $X$.
The rest of this paper develops a tool, the twisted Brauer space, to solve the problem. In
various concrete examples, the twisted Brauer space will turn out to encode interesting
geometric and arithmetic information about $X$. Perhaps the central thesis is that while
this problem would be intractable using triangulated categories, by using stable
$\infty$-categories one is able to give a precise answer which moreover accords with our
intuition: twists are classified by $1$-cocycles in automorphisms. There is a subtlety,
which is that in this setting the automorphisms really form a topological space,
and so twists are classified by $1$-cocycles in a sheaf of spaces. That this can be made
precise is a triumph of the work of Lurie, To\"en, and others on $\infty$-categories.

One might ask to classify more generally all stable $\infty$-categories $\Mrm$ such that
$\Mrm_{k^{\sep}}\we\Mod_{A\otimes_k k^{\sep}}$. An important structural theorem due to
To\"en~\cite{toen-derived} in the simplicial commutative setting and
Antieau-Gepner~\cite{ag} in the $\EE_\infty$-setting shows that these classification
problems are the same: any such $\Mrm$ is already a module category for some $k$-algebra $B$.

Note that parts of the problem of classifying \'etale twists have already been studied. For
instance, if two schemes $X$ and $Y$ become isomorphic over $k^{\sep}$, then $\Mod_Y$
is a twisted form of $\Mod_X$. Thus, the cohomology set $\Hoh^1_{\et}(\Spec k,\ShAut_X)$,
which classifies \'etale twists of $X$ as a scheme,
contributes to the answer of the problem. If two varieties $X$ and $Y$ are derived
equivalent, then \'etale twists of each of $\Mod_X$ and $\Mod_Y$ give different
interpretations for the answer.

Besides the case of schemes, another version of this problem is very well-known, although
possibly in a different guise. Suppose one attempts to find ordinary $k$-algebras $A$ such that
$A\otimes_k\overline{k}$ is Morita equivalent to $\overline{k}$. Then, every such algebra
$A$ is Morita equivalent to a central division algebra $D$ over $k$. So, the Brauer group
$\Br(k)$ classifies these algebras. This remains true in the derived world: every
$A_{\infty}$-algebra $A$ such that $A\otimes_k\overline{k}$ is derived Morita equivalent to
$\overline{k}$ is derived Morita equivalent to a central division algebra over $k$. This
result is due to To\"en~\cite{toen-derived}. The Brauer group again has a cohomological
interpretation: it is $\Hoh^2_{\et}(k,\Gm)$.

Of course, there is no reason to settle for classifying algebras over $k$. One can also
attempt to classify algebras over a scheme $X$. So, consider the problem of classifying
sheaves of quasi-coherent $A_{\infty}$algebras $A$ over $X$ such that there is an \'etale cover
$p:U\rightarrow X$ where $\StMod_{p^*A}\we\StMod_U$, where this is an equivalence of
$U$-stacks of module categories. The derived Brauer group of $X$ is obtained by taking all
such algebras and taking the quotient by derived Morita equivalence of $X$-stacks.
It turns out that the derived Brauer group is computable with cohomological methods. When
$X$ is an ordinary scheme, the derived Brauer group is
$\Hoh^2_{\et}(X,\Gm)\times\Hoh^1_{\et}(X,\ZZ)$.
If I only cared about ordinary algebras, there would be a problem at this point:
for some quasi-compact and quasi-separated schemes, not every derived Brauer
class is the class of an ordinary algebra (see~\cite{ag}*{Section 7.5}).

My point in the previous paragraph is simply that in order to obtain a cohomological
classification, which might be amenable to computation, of Azumaya algebras, it is important
to allow $A_{\infty}$-algebras.

The examples of \'etale twists of schemes and of the Brauer group show that the solution to
Problem~\ref{prob:classification} should be very interesting, and that it should be in some
way cohomological. As in the example of the previous section, it is frequently easy to
construct some examples, but showing that they are exhaustive is much more difficult, and
this is why cohomological methods are important. Such methods are already required to show, for
instance, that $\Br(\ZZ)=0$ (see~\cite{grothendieck-brauer-3}).

\subsection{The twisted Brauer space}\label{sec:tbsintro}

Let me describe the main tool of this paper in a special case.
Let $R$ be a commutative ring (or a connective commutative ring spectrum), and let $A$ be an $R$-algebra (hence,
an $A_{\infty}$-ring).

\begin{theorem}
    There is a sheaf of spaces $\ShBr^A$ on the \'etale site of $\Spec R$ with homotopy sheaves
    \begin{equation*}
        \pi_i^s\ShBr^A\iso\begin{cases}
            0                   &   \text{if $i=0$,}\\
            \ShAut_{\Mod_A}    &   \text{if $i=1$,}\\
            \THHB^0_R(A)^{\times}                 &   \text{if $i=2$,}\\
            \THHB^{2-i}_R(A)   &   \text{if $i\geq 3$,}
        \end{cases}
    \end{equation*}
    where $\THHB^*_R(A)$ is the Hochschild cohomology sheaf of $A$ over $\Spec R$.
    There is a fringed spectral sequence
    \begin{equation*}
        \Eoh_2^{p,q}=\Hoh^p_{\et}(\Spec R,\pi_q^s\ShBr^A)\Rightarrow\pi_{q-p}\ShBr^A(R),
    \end{equation*}
    which converges completely (in the sense of fringed spectral sequences) when either $A$ is a smooth
    and proper $R$-algebra or $A$ and $R$ are ordinary rings.
    The set $\Br^A(R)=\pi_0\ShBr^A(R)$ solves Problem~\ref{prob:classification}. Namely,
    every $R$-algebra $B$ such that $B$ is \'etale locally derived Morita equivalent to $A$
    determines a point of the space $\ShBr^A(R)$ and conversely. Two points $B_0$ and $B_1$ are connected by a path if and
    only if $\Mod_{B_0}\we\Mod_{B_1}$. Moreover, there is an action of the derived Brauer group
    $\Br(R)$ on $\Br^A(R)$. If $Z$ is a derived Azumaya $R$-algebra, then $[Z]\cdot[B]=[Z\otimes_R B]$.
\end{theorem}

The twisted Brauer space and the spectral sequence are generalizations of the Brauer space
and spectral sequence developed in Antieau-Gepner~\cite{ag}. Besides having a computational
tool to compute twists, the twisted Brauer space together with its action of the (untwisted)
Brauer space carries a large amount of arithmetic information. For instance, the stabilizer
of $\Mod_C$ in $\ShBr^{\PP^1}(k)$, where $C$ is a smooth projective genus $0$ curve has
enough information to determine over which fields $C$ has rational points.

When $R$ is a connective ring spectrum, the Brauer space $\ShBr(R)$ is a $2$-fold delooping
of the units spectrum $R^\times$. When $A$ is an $R$-algebra
the space $\ShBr^A(R)$ should be viewed as a $2$-fold delooping of the
spectrum of units in $\THH^*_R(A)$. Note that this is exactly the correct amount of delooping. Since $A$
is an $A_\infty$-algebra, which is the same as being an $\EE_1$-algebra, its Hochschild
cohomology $\THH^*_R(A)$ is an $\EE_2$-algebra by Deligne's conjecture, which has been
proven by many authors; see~\cite{ha}*{Section 6.1.4}. So,
the spectrum of units is a $2$-fold loopspace. The twisted Brauer space construction has the
same formal properties of $\ShBr(R)$. For example, the group
$\Aut_{\Mod_A}$ is the derived Picard group of $A$; that is, it is the group of invertible
(complexes of) $(A,A)$-bimodules. The derived Picard group was introduced by
Rickard~\cite{rickard-equivalences}, and has been studied extensively, by
Miyachi-Yekutieli~\cite{miyachi-yekutieli} and Rouquier-Zimmermann~\cite{rouquier-zimmermann}.

When $R$ is an ordinary commutative ring and $A$ is an ordinary associative $R$-algebra, or
$X$ is an ordinary $R$-scheme, then $\THH_R^{2-i}(A)=0$ (resp. $\THH_R^{2-i}(X)=0$) for
$i\geq 3$, since one can create projective (resp. locally free) resolutions.

The spectral sequence is used to show that $\Br^C(\RR)$ does indeed
have exactly three elements, to classify noncommutative
\'etale twists of curves and quadric hypersurfaces, and to classify
twists of a certain path algebra, which corresponds to noncommutative
projective space. These last twists lead to noncommutative Severi-Brauer varieties.

Let me explain briefly two of these examples.

Given an elliptic curve $E/k$, there are three interesting groups that act on $\Drm^b(E)$.
The first is the automorphism group of $E$ as a variety, which is an extension of the
automorphism group of $E$ as an elliptic curve (a finite group) by $E$ acting on itself acting via translation.
The twists by this action are homogeneous spaces for twists of $E$ as an elliptic curve.
The curve $E$ also acts on $\Drm^b(E)$
by viewing it as the moduli space of line bundles of degree $0$ over $E$. The action is then
given by tensoring with line bundles. Twists by this action lead to the twisted derived
categories $\Drm^b(E,\alpha)$ for $\alpha\in\Br(E)$. This makes sense as every such class
$\alpha$ is killed by passage to the algebraic closure of $k$.

But, there is a final group acting on $\Drm^b(E)$, which is
$\widetilde{\SL}_2(\ZZ)$, an extension of $\SL_2(\ZZ)$ by $\ZZ$.
It follows that modular representations in $\SL_2(\ZZ)$ give rise to twists of
$\Drm^b(E)$. Unlike in the other two cases, this action does not preserve the natural
$t$-structure on $\Drm^b(E)$, and hence the twists are truly derived. The interesting point
is that \emph{every} twist of $\Drm^b(E)$ is ``built out of'' four things: central simple algebras
over $k$, homogeneous spaces over twists of $E$ as an elliptic curve, the abelian
categories mentioned above, and the derived categories associated to modular representations.


The quiver $\Omega_n$ consists of two points $a$ and $b$ and $n$ arrows from $a$ to $b$.
\begin{figure}[h]
    \centering
    \begin{tikzpicture}
        \draw [fill] (-.1,.05) circle [radius=.05] node[left] {$a$};
        \draw [fill] (1.1,.05) circle [radius=.05] node[right] {$b$};
        \draw [->] (0,.1) to [out=45, in=135] node[above] {$s_1$} (1,.1);
        \draw [->] (0,.05) to [out=30, in=150] (1,.05);
        \draw [->] (0,0) to [out=-30, in=210] node[below] {$s_n$} (1,0) ;
        \draw [dotted] (.5,.17) to (.5,-.13);
    \end{tikzpicture}
    \caption{The quiver $\Omega_n$.}
\end{figure}
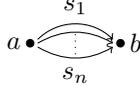
Kontsevich and Rosenberg showed that the path algebra $k\Omega_n$ is derived equivalent to
the derived category of coherent sheaves on noncommutative projective space $\NP^{n-1}$.
For $n\geq 3$, $\NP^{n-1}$ and $\PP^{n-1}$ are not derived equivalent, so these spaces are
new from the perspective of noncommutative algebraic geometry above.
However, Miyachi and Yekutieli~\cite{miyachi-yekutieli}*{Corollary 0.4} computed the automorphisms of the derived category
of $k\Omega_n$, showing that it is an extension of $\ShPGL_n(k)$. Using their calculation, the
work below shows that there is one twist of $\Drm^b(k\Omega_n)$ for each classical
Severi-Brauer variety over $k$. Thus, the twists of $k\Omega_n$ are
noncommutative Severi-Brauer varieties.

By restricting attention to simplicial commutative rings, as for instance
used by To\"en~\cite{toen-derived} and To\"en-Vaqui\'e~\cite{toen-vaquie}, it is possible to
use the fppf topology instead of the \'etale topology. The theory below carries over without
change to the simplicial setting.

In Section~\ref{sec:tbs}, the necessary background is reviewed and the definition and first
properties of the twisted Brauer space are studied. The spectral
sequence that computes the homotopy of the twisted Brauer space is constructed in
Section~\ref{sec:ss}. This is used
to give a complete description of $\Br^{\PP^1}(\RR)$. In
Section~\ref{sec:moritadescent}, the problem of when it is enough to check
derived Morita equivalence over $\overline{k}$ is considered.
Several examples are studied in Section~\ref{sec:geoexamples}.

\subsection{Acknowledgments}

This paper would not exist without my collaboration with David Gepner, who I would
like to thank for patiently explaining to me many things about
$\infty$-categories during the writing of~\cite{ag}. His perspective on higher algebra is
present everywhere in this work. I also thank Rapha\"el Rouquier for several useful conversations.

\section{The twisted Brauer space}\label{sec:tbs}

\subsection{Derived Morita theory}

Recall from the introduction that if $A$ is an $A_{\infty}$-algebra, then $\Mod_A$ denotes the
stable $\infty$-category of right $A$-modules. This is a large $\infty$-category because
it is complete and cocomplete. The subcategory $\Mod_A^{\ccc}$ is the small stable
$\infty$-category of compact $A$-modules. For a scheme $X$, $\Mod_X$ denotes the stable
$\infty$-category of complexes of $\Oscr_X$-modules with quasi-coherent cohomology sheaves.
Then, $\Mod_X^{\ccc}$, the subcategory of compact objects, is the same as the
$\infty$-category of perfect complexes on $X$, at least when $X$ is quasi-compact and
quasi-separated. See~\cite{bondal-vandenbergh}.

The $\infty$-categories $\Mod_A$ and $\Mod_X$, besides being stable, are also presentable
$\infty$-categories, which is equivalent to saying that their homotopy categories are have
all coproducts, are locally small, and are $\kappa$-compactly generated for some regular
cardinal $\kappa$ (see~\cite{ha}*{Corollary 1.4.4.2}. This fact follows from Lurie~\cite{ha} in the case of $\Mod_A$, and
from~\cite{bondal-vandenbergh} when $X$ is quasi-compact and quasi-separated, because in
that case $\Mod_X\we\Mod_A$ for an appropriate choice of $A$. Presentability ensures that
$\Mod_A$ and $\Mod_X$ have all small limits and colimits and that they
can be described by a set of generators in a reasonable way. For details,
see~\cite{htt}*{Chapter 5}. If $A$ is an $R$-algebra, where $R$ is a commutative ring
spectrum, then $\Mod_A$ is enriched over $R$, in the sense that the mapping spectra in $\Mod_A$ are
naturally $R$-modules. By an $R$-linear category, I will mean a stable presentable
$\infty$-category enriched over $\Mod_R$. There is a symmetric monoidal $\infty$-category $\Cat_R$ whose objects are $R$-linear
categories and whose morphisms are $R$-linear functors that have right adjoints.

There are two points of derived Morita theory to bear in mind for the paper below. First, if
$A$ and $B$ are $R$-algebras, then
any $R$-linear functor $F:\Mod_A\rightarrow\Mod_B$ in $\Cat_R$ is determined by an
$A^{\op}\otimes B$-module $F(A)$. Moreover, there are natural equivalences
\begin{equation*}
    \Fun_R^{\Lrm}(\Mod_A,\Mod_B)\we\Mod_{A^{\op}}\otimes_{\Mod_R}\Mod_B\we\Mod_{A^{\op}\otimes_R B},
\end{equation*}
where $\Fun_R^{\Lrm}(-,-)$ denotes the functor $\infty$-category of left adjoint $R$-linear
functors.

The second point is that if $E$ is a compact generator of any $R$-linear stable
$\infty$-category $\Mrm$, then the mapping spectrum out of $E$ induces an equivalence
\begin{equation*}
    \Map_{\Mrm}(E,-):\Mod_A\rightarrow\Mod_{\End_R(E)^{\op}}
\end{equation*}
by Schwede-Shipley~\cite{schwede-shipley}.
The converse is also true. In particular, the result of Bondal and van den Bergh
says that there is a compact generator of $\Mod_X$ when $X$ is quasi-compact and
quasi-separated, so $\Mod_X\we\Mod_A$ for some $A_{\infty}$-algebra $A$.

These results should be compared in two directions to more familiar facts. First, they are
essentially a translation into the world of stable $\infty$-categories of facts that are
true for abelian categories of modules, from which the appellation Morita originated.
Second, for a scheme, functors $\Mod_X\rightarrow\Mod_Y$ are determined by complexes on
$X\times Y$. For fully faithful functors, Orlov proved this result for functors of the
derived categories $\Drm^b(X)\rightarrow\Drm^b(Y)$; it is a very difficult and important
theorem. At the level of $\infty$-categorical models, it is due to Ben-Zvi, Francis, and
Nadler~\cite{bzfn}, while To\"en proves it for dg models~\cite{toen-morita}.

As a last point in this section of background, if $\Mrm$ is an $R$-linear category, and if
$S$ is a commutative $R$-algebra, then one can base-change $\Mrm$ up to $S$ via
$\Mrm_S=\Mod_S\otimes_{\Mod_R}\Mrm$.

\subsection{The definition}

Let $R$ be a connective commutative ring spectrum, and let $\Shv_R^{\et}$ be the big \'etale topos over
$\Spec R$. If $S$ is a connective commutative $R$-algebra, an $S$-linear category $\Mrm$ is
said to satisfy \'etale hyperdescent if
for every connective commutative $S$-algebra $T$ and every \'etale hypercover $T\rightarrow
U^\bullet$ of $T$, the induced morphism
\begin{equation*}
    \Mrm_T\rightarrow\lim_{\Delta}\Mrm_{U^\bullet}
\end{equation*}
is an equivalence.
There is a stack of large $\infty$-categories $\StCat^{\desc}$ over $\Spec R$ that
classifies linear categories with \'etale hyperdescent and left adjoint functors
between them~\cite{dag11}*{Theorem 7.5}. Write $\ShPr$ for the underlying sheaf of spaces.
For details, see Antieau-Gepner~\cite{ag}*{Section 6}. 

Suppose now that $Z$ is in $\Shv_R^{\et}$, and let $\alpha:Z\rightarrow\ShPr$ be a map of
sheaves. The corresponding linear category with descent, or,
equivalently, stack of linear categories, will be denoted $\StMod^\alpha$. The
$\infty$-category of sections over $f:\Spec S\rightarrow Z$ is the $S$-linear category
$\Mod_S^{f\circ\alpha}$ classified by $f\circ\alpha$ by Yoneda's lemma. By definition, the
$\infty$-category $\Mod_X^{\alpha}$ of sections over a sheaf $X$ over $Z$ is
\begin{equation*}
    \Mod_X^{\alpha}=\lim_{f:\Spec S\rightarrow X}\Mod_S^{f\circ\alpha}.
\end{equation*}
For instance, let $\Oscr:\Spec R\rightarrow\ShPr$ send $\Spec S$ to $\Mod_S$. Then,
$\Mod_X^\Oscr$ is the stable $\infty$-category of quasi-coherent $\Oscr_X$-modules.
The properties of this construction of sheaves have been studied extensively in~\cite{bzfn},
\cite{dag11}, and \cite{ag}.

For an object $f:X\rightarrow Z$ of
$\Shv_Z^{\et}=\left(\Shv_R^{\et}\right)_{/Z}$, there is a pullback stack $f^*\alpha$. Say
that a stack of linear categories $\beta:X\rightarrow\ShPr$ over $X$ is \'etale locally
equivalent to $f^*\alpha$ if there is an \'etale cover $p:U\rightarrow X$ such that
$p^*\beta\we p^*f^*\alpha$ as stacks of linear categories over $U$. There is a
subspace $\ShBr^\alpha(X)$ of $\ShPr(X)$ of stacks of linear categories that are \'etale locally
equivalent to $f^*\alpha$.

\begin{lemma}
    The presheaf $\ShBr^\alpha$ on $\Shv_Z^{\et}$ is an \'etale sheaf.
    \begin{proof}
        The presheaf is the same as the sheafification of the point $\alpha$ in $\ShPr|_Z$.
    \end{proof}
\end{lemma}

\begin{definition}
    The sheaf of spaces $\ShBr^{\alpha}$ is called the $\alpha$-twisted Brauer sheaf. For
    a sheaf $X$, $\ShBr^{\alpha}(X)$ is the $\alpha$-twisted Brauer space of $X$. The
    pointed set $\pi_0\ShBr^{\alpha}(X)$ is the $\alpha$-twisted Brauer set of $X$, and it
    will be written $\Br^{\alpha}(X)$ in the sequel.
\end{definition}

To summarize in a fast and loose way in a familiar setting, if $X$ is a $k$-variety, where
$k$ is a field, and if
$A$ is an ordinary associative $k$-algebra, then the twisted Brauer set $\Br^A(X)$ classifies
sheaves quasi-coherent dg algebras $\Bscr$ that are \'etale locally derived Morita equivalent on $X$ to
$\Oscr_X\otimes_k A$. This is fast because it has yet to be observed that elements
of $\Br^A(X)$ actually correspond to algebras, although this is true; see the next
section. The only looseness in this description is that the \'etale-local
Morita equivalence is an equivalence of the \emph{stacks} of modules. See
Remark~\ref{rem:stackymorita} at the end of the section.

For example, if $\Oscr$ classifies the stack of quasi-coherent modules over $Z$, then
$\ShBr^\Oscr=\ShBr$, the Brauer sheaf studied in~\cite{ag}.


\begin{example}
    Suppose that $A$ is an associative $S$-algebra. Then, the stack $\StMod^A$ is the stack of
    linear categories whose $\infty$-category of sections over a connective commutative
    $S$-algebra $T$ is $\Mod_{T\otimes_S A}$. In this case, the twisted Brauer sheaf is
    $\ShBr^A$.
\end{example}

\begin{example}
    Suppose that $X$ is a scheme over $\Spec S$. Then, $\StMod^X$ is the stack of linear
    categories over $\Spec S$ whose category of sections over $T$ is $\Mod_{X_T}$, where
    $X_T=X\times_{\Spec S}\Spec T$. Here $\Mod_{X_T}$ is the stable $T$-linear
    $\infty$-category with homotopy category equivalent to $\Drm_{\qc}(X_T)$, the derived
    category of complexes of $\Oscr_{X_T}$-modules with quasi-coherent cohomology.
    This slightly unusual notation is meant to emphasize that $\StMod^X$ is viewed not as a stack over
    $X$ but over $\Spec S$. 
    The associated twisted Brauer sheaf is $\ShBr^X$.
    Note that if $X$ is quasi-compact and quasi-separated,
    then by the results of~\cite{bondal-vandenbergh}, this is a special case of the previous
    example. When $X\rightarrow\Spec S$ is smooth, then the elements of $\ShBr^X$
    may viewed as \'etale twists of $\Drm^b(X)$. In general, they should be viewed as either
    twists of $\Mod_X$ or $\Perf_X$.
\end{example}

It is not clear that, in general, $\ShBr^\alpha$ is a sheaf of small spaces. However, in
most cases of interest, and all cases considered in this paper, it is. To prove this, we
need a lemma first, which will be of use later in the paper for computing twisted
Brauer spaces.

\begin{lemma}\label{lem:autoeq}
    The sheaf $\ShBr^\alpha$ is equivalent to the classifying sheaf of the sheaf of
    autoequivalences of the stack $\alpha$.
    \begin{proof}
        By definition, any two points of $\ShBr^\alpha(X)$ are \'etale locally connected. It
        follows that the homotopy sheaf $\pi_0^s\ShBr^\alpha$ is just a point. There is
        an obvious morphism $\ShB\Shaut(\alpha)\rightarrow\ShBr^\alpha$. So, it suffices to
        compute the homotopy sheaves of the loopspace $\Omega\ShBr^\alpha$ at the point
        $\alpha$. But, these are just the equivalences from the stack $\alpha$ to $\alpha$, as
        desired.
    \end{proof}
\end{lemma}

We say that $\alpha:Z\rightarrow\ShPr$ classifies a stack of compactly generated linear
categories if $\Mod_S^\alpha$ is compactly generated for every $\Spec S\rightarrow Z$ and
every connective commutative $R$-algebra $S$. Note that this hypothesis does not imply that,
for instance, $\Mod_Z^\alpha$ is compactly generated. However, if $Z$ is a quasi-compact and
quasi-separated derived scheme, then the methods of Lurie~\cite{dag11}*{Section 6}
can be used to show that $\Mod_Z^\alpha$ is compactly generated.

\begin{proposition}
    Suppose that $\alpha$ classifies a stack of compactly generated linear categories over a
    sheaf $Z$. Then, $\ShBr^\alpha$ is a sheaf of small spaces.
    \begin{proof}
        By the previous lemma, it is enough check that $\Shaut(\alpha)$ is a sheaf
        of small spaces, which we can check on affines $\Spec S\rightarrow Z$, and by
        hypothesis $\Mod_S^\alpha$ is compactly generated.
        The space of sections over $\Spec S\rightarrow Z$ is
        a subspace of the $\infty$-category of functors
        $\Mod_S^\alpha\rightarrow\Mod_S^\alpha$
        that preserve the subcategory of compact objects $\Mod_S^{\alpha,\ccc}$. Thus, it is a subspace of functors
        $\Mod_S^{\alpha,\ccc}\rightarrow\Mod_S^{\alpha,\ccc}$, where
        $\Mod_S^{\alpha,\ccc}$ is the full subcategory of $\Mod_S^\alpha$ of compact
        objects. Write $-\alpha$ for the
        opposite linear category. Then, $\Shaut(\alpha)(S)$ is a subspace of
        $\Mod_S^{-\alpha,\ccc}\otimes\Mod_S^{\alpha,\ccc}$, which is a small idempotent
        complete stable $\infty$-category. It follows that $\Shaut(\alpha)$ is a sheaf of
        small spaces.
    \end{proof}
\end{proposition}

\begin{remark}\label{rem:stackymorita}
    It is worth noting that there is a subtlety required when defining the Brauer group via
    Morita equivalences. If $X$ is a scheme with an automorphism $\sigma$ such that
    $\sigma^*\alpha\neq\alpha$ for some Brauer class $\alpha$, then it is clear that the
    categories of $\alpha$-twisted coherent sheaves
    $\Mod_X^\alpha$ and $\Mod_X^{\sigma*\alpha}$ are equivalent. So, the Brauer group is a
    finer invariant than this sort of derived equivalence. I learned of this issue from
    C\u{a}ld\u{a}raru's thesis~\cite{caldararu}*{Example 1.3.16}. However, the underlying stacks
    $\StMod^{\alpha}$ and $\StMod^{\sigma^*\alpha}$ are \emph{not} equivalent over $X$. The Brauer
    group classifies Azumaya algebras up to the Morita equivalence classes of the stacks of
    their modules. This perspective is implicit in~\cite{toen-derived} and~\cite{ag}.
    The hypothesis
    $\Mod_X^\alpha\we\Mod_X^{\sigma^*\alpha}$ is rather strange from the perspective of
    stacks: it is like saying that two coherent sheaves have isomorphic vector spaces of global
    sections. The correct definition of stacky Morita equivalence is built into
    $\ShBr^\alpha(X)$.
\end{remark}

\subsection{Algebras and the twisted Brauer space}

If $\alpha:Z\rightarrow\ShPr$ is a stack of linear categories, and if $f:X\rightarrow Z$ is a
map of sheaves, then an object $x$ of $\Mod_X^\alpha$ is perfect if its restriction to
$\Mod_T^\alpha$ for every $\Spec T\rightarrow X$ compact.
A set of perfect objects $\Gamma$ perfectly generates $\Mod_X^\alpha$ if its restriction to every
affine is a set of compact generators. It globally generates if it perfectly generates, if
the objects are compact, and if it generates $\Mod_X^\alpha$.

Lurie has shown in~\cite{dag11}*{Theorem 6.1} that if $\Mod^\alpha_S$ is an $S$-linear
category with descent that is \'etale locally compactly
generated, then $\Mod_S^\alpha$ is globally generated. In~\cite{ag}*{Theorem 6.17}, Gepner and I showed that if
$\Mod^\alpha_S$ is \'etale locally compactly generated by a \emph{single} compact object, then
$\Mod_S^\alpha$ is globally generated by a single compact object. To\"en proved similar
theorems for simplicial commutative rings in~\cite{toen-derived}, although with somewhat
different methods.

\begin{proposition}
    Suppose that $\alpha:Z\rightarrow\ShPr$ classifies a stack of linear categories, and fix
    a morphism $f:X\rightarrow Z$, where $X$ is a quasi-compact and quasi-separated derived
    scheme. Then, if $\Mod_X^{\alpha}$ is globally generated, so is $\Mod_X^\beta$ for
    every $\beta\in\ShBr^\alpha(X)$.
    Similarly, if $\Mod_X^{\alpha}$ is globally generated by a single object, then so is
    $\Mod_X^\beta$ for every $\beta\in\ShBr^\alpha(X)$.
    \begin{proof}
        These statements follow from Lurie~\cite{dag11}*{Theorem 6.1}
        and Antieau-Gepner~\cite{ag}*{Theorem~6.17}, respectively. For instance, if
        $\Mod_X^\alpha$ is globally generated by a single object, then, using the \'etale local equivalence of
        $\alpha$ and $\beta$ over $X$, it follows that $\Mod_X^{\beta}$ is \'etale locally
        compactly generated by a single compact object. Now, apply~\cite{ag}*{Theorem 6.17}.
    \end{proof}
\end{proposition}

As a corollary, in the case of generation by a single object, the stacks $\StMod^\beta$ are
stacks of modules for a quasi-coherent algebra. This is a useful thing to know, as it can
make it easier to compute Hochschild cohomology and other invariants of the categories.

\begin{corollary}
    Let $\Ascr:Z\rightarrow\ShPr$ classify the stack of $\Ascr$-modules for a quasi-coherent
    algebra $\Ascr$ over $Z$. If $p:X\rightarrow Z$ is a quasi-compact and quasi-separated derived
    scheme, then for every $\beta\in\ShBr^\Ascr(X)$, the stack $\StMod^{\beta}$ is equivalent to the stack of modules
    $\StMod^\Bscr$ for a quasi-coherent $\Oscr_X$-algebra $\Bscr$.
\end{corollary}

The corollary is a twisted and derived form of the $\Br=\Br'$ question of Grothendieck, which
asks if every cohomological Brauer class comes from an Azumaya algebra. This turns out to be
false if one only considers ordinary Azumaya algebras over an ordinary scheme. It is
necessary in
certain cases to allow derived Azumaya algebras as well. On the other hand, if $X$ has an
ample line bundle, then a theorem of Gabber (see de Jong~\cite{dejong}) showed that
$\Br(X)=\Hoh^2_{\et}(X,\Gm)_{\tors}$. If in addition $X$ is regular and noetherian, then the computations
of~\cite{ag}*{Section 7} show that every derived Azumaya algebra on $X$ is Morita equivalent
to an ordinary Azumaya algebra.

\begin{question}\label{q:grothendieck}
    Suppose that $R$ is a regular, noetherian ring, and let $A$ be an ordinary associative
    $R$-algebra. Is every element $\beta\in\ShBr^A(\Spec R)$ derived Morita equivalent to an
    ordinary $R$-algebra $B$?
\end{question}

Although this seems like a difficult question in general, this paper gives a positive answer
for quadric hypersurfaces and noncommutative projective spaces.

\subsection{The action of the Brauer group and $\alpha$-twisted sheaves}

There is an action of the Brauer space $\ShBr$ on $\ShBr^\alpha$ for any $\alpha$. Indeed,
if $\beta:X\rightarrow\ShPr$ is \'etale locally equivalent to $\alpha$, and if
$\gamma:X\rightarrow\ShPr$ is \'etale locally equivalent to $\Oscr:X\rightarrow\ShPr$, then
the tensor product $\gamma\otimes\beta$ is \'etale locally equivalent to $\alpha$, since, if
$\Spec S\rightarrow X$ is a map from an affine on which $\beta$ is equivalent to $\alpha$
and $\gamma$ is equivalent to $\Oscr$, one sees that
$$\Mod_S^{\gamma}\otimes_{\Mod_S}\Mod_S^{\beta}\we\Mod_S\otimes_{\Mod_S}\Mod_S^\alpha\we\Mod_S^\alpha.$$

A special case of this action has already gained a great deal of attention under a different guise, namely as
derived categories of twisted sheaves (see~\cite{caldararu} or~\cite{lieblich}). Suppose that $X$ is an ordinary scheme and that
$\alpha\in\Br'(X)=\Hoh^2(X,\Gm)_{\tors}$. One can represent $\alpha$ as a $2$-cocycle
$(\alpha_{ijk})$ over
some \'etale cover $\{U_i\}_{i\in I}$ of $X$. An $\alpha$-twisted coherent sheaf consists of
a coherent $\Oscr_{U_i}$-module $\Fscr_i$ for each $i$ and an isomorphism
$\theta_{ij}:\Fscr_i|_{U_{ij}}\rightarrow\Fscr_j|_{U_{ij}}$ such that
\begin{equation*}
    \theta_{ki}\circ\theta_{jk}\circ\theta_{ij}
\end{equation*}
is multiplication by $\alpha_{ijk}$ on $\Fscr_i|_{U_{ijk}}$. The $\alpha$-twisted coherent
sheaves form an abelian category and so one can speak of complexes of $\alpha$-twisted
coherent sheaves and obtain a derived category $\Drm^b(X,\alpha)$.

More generally, one can consider the stable $\infty$-category $\Mod_X^{\alpha}$ of complexes
of $\alpha$-twisted $\Oscr_X$-modules with quasi-coherent ($\alpha$-twisted) cohomology
sheaves. The subcategory of compact objects may be identified with the category of complexes
of $\alpha$-twisted perfect complexes.

If $\alpha\in\Br'(X)$, write $\StMod^{X,\alpha}$ for the stack (over
$\Spec R$) of $\alpha$-twisted sheaves on $X$.
If $X$ is regular and noetherian,
$\Drm^b(X,\alpha)\we\Ho(\Mod^{X,\alpha,\perf}_R)$, the homotopy category of the stable
$\infty$-category of $\alpha$-twisted perfect complexes on $X$.
These are the compact objects in $\Mod^{X,\alpha}$.

Let $X$ be a scheme over an ordinary commutative ring $S$.
For an arbitrary $\alpha\in\Br'(X)$, it will not be the case that $\StMod^{X,\alpha}$ will be
\'etale locally equivalent to $\StMod^X$ over $\Spec S$. For instance, if $X$ is a K3
surface over an algebraically closed field $k$, then $\Drm^b(X)$ is not equivalent to
$\Drm^b(X,\alpha)$ for any $\alpha\neq 0$ in $\Br'(X)$. The stacks of linear categories
$\StMod^X$ and $\StMod^{X,\alpha}$ are in general only \'etale locally equivalent on $X$.
However, if $\alpha\in\Br(S)$, we can pull-back via the structure morphism
$p:X\rightarrow\Spec S$ to obtain $p^*\alpha$. Then, $\StMod^{X,p^*\alpha}$ is \'etale locally
Morita equivalent to $\StMod^X$ over $\Spec S$.

\begin{proposition}
    The action of $\alpha\in\Br(S)$ on $\Br^X(S)$ sends $\StMod^X$ to $\StMod^{X,p^*\alpha}$.
    \begin{proof}
        Write $\pi:X\rightarrow\Spec S$.
        By definition, $\alpha\cdot\StMod^X$ is the stack that sends $f:\Spec T\rightarrow\Spec
        S$ to the $T$-linear category
        \begin{equation*}
            \Mod^{\alpha}_T\otimes_{\Mod_T}\Mod^{X}_T\we\Mod^{\alpha}_T\otimes_{\Mod_T}\Mod_{X_T}.
        \end{equation*}
        On the other hand, $\StMod^{X,p^*\alpha}$ is the stack that sends $f$ to the
        $T$-linear category $\Mod_{X_T}^{p^*\alpha}$. There is a natural map from
        $\Mod^{\alpha}_T\otimes_{\Mod_T}\Mod_{X_T}$ to $\Mod_{X_T}^{p^*\alpha}$, which is an
        equivalence \'etale locally on $\Spec T$. It follows that it is already an
        equivalence by descent. Taking $T=S$, the proposition follows.
    \end{proof}
\end{proposition}

The action will be given a cohomological interpretation at the end of
Section~\ref{subsec:ss}.

\begin{corollary}\label{cor:stabilizers}
    If the pullback $p^*\alpha$ is zero in $\Br(X)$, then $\alpha$ stabilizes $\StMod^X$.
\end{corollary}

I conjecture that the converse is true. The conjecture will be verified in various cases
throughout the paper, including for smooth projective varieties $X$ over a field with
$\omega_X$ ample or anti-ample.

\begin{conjecture}[Stabilizer conjecture]\label{conj:stabilizers}
    If $\alpha\in\Br(S)$ stabilizes $\StMod^X$, then
    $\alpha\in\ker(\Br(S)\rightarrow\Br(X))$.
\end{conjecture}

To conclude the section, I include a more formal structural remark.
If $A$, $B$, $C$, and $D$ are $S$-algebras, and if $C$ is \'etale locally Morita equivalent to $A$ and $D$ is
\'etale locally Morita equivalent to $B$, then $C\otimes_S D$ is \'etale locally Morita equivalent to
$A\otimes_S D$. Thus, there are natural products
$\ShBr(-;A)\times\ShBr(-;B)\rightarrow\ShBr(-;A\otimes_S B)$ of sheaves of spaces over
$\Spec S$. The Brauer sheaf $\ShBr$ is an $\EE_\infty$-algebra object in
$\Shv_S^{\et}$ by~\cite{ag}*{Corollary 7.5}, which means that it is a sheaf of group-like $\EE_\infty$-spaces.

\begin{proposition}
    If $A$ is an $S$-algebra, then $\ShBr^A$ is a module for the $\EE_\infty$-algebra
    object $\ShBr$, and thus can be viewed as an element of $\Mod_{\ShBr}(\Shv_S^{\et})$.
    There is a natural equivalence
    \begin{equation*}
        \ShBr^A\otimes_{\ShBr}\ShBr^B\we\ShBr^{A\otimes_S B}.
    \end{equation*}
    \begin{proof}
        The claim that $\ShBr^A$ is a module for $\ShBr$ follows from the symmetric
        monoidal structure on the sheaf $\ShPr$. To prove the second claim, it is enough to
        note that if $T$ is a connective commutative $S$-algebra and $C$ is a $T$-algebra
        that is \'etale locally equivalent to $A\otimes_S B\otimes_S T$, then $C$ is \'etale
        locally equivalent to the tensor product of an algebra in $\ShBr^A(T)$ and an
        algebra in $\ShBr^B(T)$.
    \end{proof}
\end{proposition}

\section{The descent spectral sequence}\label{sec:ss}

\subsection{Fringed spectral sequences}\label{sec:fringed}

To consider carefully what happens in the fringed spectral sequences that appear when doing
descent spectral sequences, it is useful to first consider the long exact sequence of
homotopy groups associated to a fibration $p:X\rightarrow Y$ of pointed spaces. Let $b\in Y$
be the basepoint, and let $f\in F$, where $F=p^{-1}\{b\}$. Then, there is a sequence of
homotopy groups and pointed homotopy sets
\begin{equation*}
    \rightarrow\cdots\pi_2(Y,b)\rightarrow\pi_1(F,f)\rightarrow\pi_1(X,f)\rightarrow\pi_1(Y,b)\rightarrow\pi_0(F,f)\rightarrow\pi_0(X,f)\rightarrow\pi_0(Y,b),
\end{equation*}
where $\pi_0(-,f)$ is the set of path
components pointed by $f$. This sequence is exact in the following sense:
\begin{itemize}
    \item   at any place $\pi_i(F,f)$, $\pi_i(X,f)$, or $\pi_i(Y,b)$ where $i>0$, it is
        exact in the usual sense that $\ker=\im$;
    \item   the image of $\pi_2(Y,b)$ is in the center of $\pi_1(F,f)$;
    \item   there is an action of $\pi_1(Y,b)$ on $\pi_0(F,f)$ such that two elements of
        $\pi_0(F,f)$ agree in $\pi_0(X,f)$ if and only if they are in the same orbit;
    \item   the map $\pi_1(Y,b)\rightarrow\pi_0(F,f)$ induces a bijection between
        $\pi_1(Y,f)/\pi_1(X,f)$ and the orbit of the point $f$ in $\pi_0(F,f)$;
    \item   a point $g\in\pi_0(X,f)$ goes to $b$ in $\pi_0(Y,b)$ if and only if it is in the
        image of $\pi_0(F,f)\rightarrow\pi_0(X,f)$.
\end{itemize}
The main information that this
sequence does not see is the fact that the fibers of $\pi_0(X,f)\rightarrow\pi_0(Y,b)$ can
vary widely over different points of $\pi_0(Y,b)$ and can be empty, so that in particular
$\pi_0(X,f)\rightarrow\pi_0(Y,b)$ might not be surjective.

Now, let $\cdots\rightarrow X_n\rightarrow X_{n-1}\rightarrow\cdots X_0\rightarrow *$ be a
sequence of fibrations of pointed spaces, where $X_n$ is pointed by $f_n$. Let $f$ be the
point of $X=\lim X_n$ that is the inverse limit of the points $f_n$. Write $F_n$ for the
homotopy fiber of $X_n\rightarrow X_{n-1}$ over $f_{n-1}$. Bousfield and Kan~\cite{bousfield-kan}*{Section
IX.4} created a spectral sequence that converges
conditionally to $\pi_*X$ by rolling up all of the fibration sequences
$F_n\rightarrow X_n\rightarrow X_{n-1}$ into a generalized triple, generalized in the sense
that some terms are not abelian groups. Without going into many details, there is a fringed spectral
sequence
\begin{equation*}
    \Eoh_2^{s,t}=\pi_{t-s}(F_t,f_t)\Rightarrow\pi_{t-s}X,
\end{equation*}
fringed in the sense that $\Eoh_r^{s,s}$ is just a pointed set, and $\Eoh_r^{s,s+1}$ is a
possibly non-abelian group.

The differential $d_r$ in the $\Eoh_r$-page of this spectral sequence has degree
$(r,r-1)$. (Note that Bousfield and Kan index
the spectral sequence differently, beginning instead with $\Eoh_1^{s,t}=\pi_{t-s}F_s$.)
When $t-s>0$, the $\Eoh_{r+1}^{s,t}$ term is computed in the usual way from
$\Eoh_r^{s,t}$, as cycles modulo boundaries. When $t-s=0$, there is not only a differential
with target $\Eoh_r^{s,t}$, but the source $\Eoh_r^{s-r,t-r+1}$ acts on $\Eoh_r^{s,t}$, and
$\Eoh_{r+1}^{s,t}$ is the orbit space of this action. The meaning of convergence is clear
when $t-s>0$. When $t-s=0$, there is a filtration of $\pi_0 X$ as a pointed set. This
means that there is a sequence of inclusions of pointed sets
$$*\subseteq\cdots\subseteq F_{s+1}\pi_0 X\subseteq F_s\pi_0X\subseteq\cdots\subseteq
F_0\pi_0X=\pi_0X$$
and the successive quotients $F_s\pi_0X/F_{s+1}\pi_0 X$ are bijective to
$\Eoh_{\infty}^{s,s}$ as pointed sets. The filtration on $\pi_iX$ has the same indexing.
Namely, there is a decreasing filtration $F_s\pi_iX$ and
$F_s\pi_iX/F_{s+1}\pi_iX\iso\Eoh_{\infty}^{s,s+i}$, when the spectral sequence converges.

The reader is warned that the convergence of this spectral sequence is in general
only conditional. However, the spectral sequence will converge completely in all cases considered in
this paper (for $t-s>0$).
For a discussion of convergence of these spectral sequence,
see~\cite{bousfield-kan}*{Section IX.5}. It makes sense only for the terms abutting to
$\pi_iX$ where $i>0$, where it coincides with the usual notion of convergence.
In general, more work is needed to get a handle on $\pi_0X$, which is
the case of greatest interest in this paper. However, the complete convergence for $t-s>0$
will often give crucial information for understanding what happens for $\pi_0X$.

The descent spectral sequence, sometimes called the Brown-Gersten spectral
sequence~\cite{brown-gersten}, associated to a sheaf of spaces on a topos is a special
case of the spectral sequence associated to a tower of fibrations. Let $F$ be a sheaf of
pointed spaces, which is to say an object of $\Shv_R^{\et}$. The construction below works
for any object in any $\infty$-topos. However, convergence is a more delicate question,
closely related to notion of hypercompleteness discussed in~\cite{htt}*{Section 6.5}.
The Postnikov tower of $F$ as a
sheaf is obtained via the truncations of $F$
\begin{equation*}
    F\rightarrow\cdots\tau_{\leq n}F\rightarrow\tau_{\leq n-1}\rightarrow\cdots\tau_{\leq
    0}F\rightarrow *,
\end{equation*}
and the fiber of $\tau_{\leq n}F\rightarrow\tau_{\leq n-1}F$ is the Eilenberg-MacLane sheaf $K(\pi_n^sF,n)$,
which has homotopy sheaves $\pi_n^sF$ in degree $n$ and $0$ (or a point) elsewhere. Let $X$
be another sheaf in $\Shv_R^{\et}$. Then, the sequence
\begin{equation*}
    \cdots\rightarrow\map(X,\tau_{\leq n}F)\rightarrow\map(X,\tau_{\leq
    n-1}F)\rightarrow\cdots,
\end{equation*}
is a tower of fibrations which, in good cases, and all cases in this paper, has inverse
limit the space $\map(X,F)$. The spectral associated to this tower has the form
\begin{equation*}
    \Eoh_2^{p,q}=\pi_{q-p}\map(X,K(\pi_q^sF,q))\Rightarrow\pi_{q-p}\map(X,F).
\end{equation*}
Since $K(\pi_q^sF,q)$ is an infinite loop space (at least for $q>1$),
$$\pi_{q-p}\map(X,K(\pi_q^sF,q))\we\pi_0\map(X,K(\pi_q^sF,p)).$$

Suppose for a moment that $C$ is a small category with a Grothendieck topology and that $X$
is an object of $\Shv_{\Nrm C}$ and $A$ is a sheaf of abelian groups on $C$. Then,
Lurie shows~\cite{htt}*{Remark 7.2.2.17} that
\begin{equation*}
    \pi_0\map(X,K(A,n))\we\Hoh^n(X,A),
\end{equation*}
where $\Hoh^n(X,A)$ denotes the usual cohomology group of $X$ with coefficients in $A$.
Since the small \'etale site over a connective commutative $R$-algebra $S$ is equivalent to
the nerve of the small \'etale site over $\pi_0 S$, it follows that
if $X=\Spec S$, then the groups
\begin{equation*}
    \pi_{q-p}\map(\Spec S,K(\pi_q^sF,q))\we\pi_0\map(\Spec
    S,K(\pi_q^sF,p))\we\Hoh^p_{\et}(\Spec\pi_0S,\pi_q^sF).
\end{equation*}
This has the following generalization to schemes.

\begin{proposition}
    Let $R$ be an ordinary commutative ring, and let $X$ be an ordinary $R$-scheme, viewed
    as an object of the $\infty$-topos $\Shv_{\Hrm R}^{\et}$. If $A$ is an abelian
    group object in the underlying discrete topos, then
    \begin{equation*}
        \pi_0\map(X,K(A,n))\iso\Hoh^n_{\et}(X,A),
    \end{equation*}
    where $\Hoh^n_{\et}(X,A)$ denotes the usual \'etale cohomology group of $X$ with
    coefficients in $A$.
    \begin{proof}
        The Eilenberg-MacLane sheaf is hypercomplete, so one can compute the group
        $\pi_0\map(X,K(A,n))$ with a suitably nice \'etale hypercover of $X$ that will
        also compute the group $\Hoh^n_{\et}(X,A)$. Assuming that
        this hypercover consists of disjoint unions of affine schemes, the observation above
        that the statement is true for affine schemes shows that the proposition is true by
        comparing the \v{C}ech complexes.
    \end{proof}
\end{proposition}

\subsection{The spectral sequence}\label{subsec:ss}

Here and in the rest of the paper I will abuse notation and use $\Mod_A$ and $\StMod^A$
interchangeably. There is no danger of confusion or error, since $\Mod_A$ is an $R$-linear
category with descent, so that $\StMod^A$ can be constructed from $\Mod_A$, and vice versa.

In this paper, the main objects of interest are $\ShBr^X$, where $X$ is a smooth proper
variety over an ordinary commutative ring $R$, which is a special case of $\ShBr^A$ where
$A$ is a smooth and proper $R$-algebra. The strategy for actually computing $\ShBr^X(R)$ is
to determine the sheaf of spaces $\Shaut_{\Mod_X}$, use the fact that $\ShBr^X$ is the classifying
sheaf of $\Shaut_{\Mod_X}$, and use the descent spectral sequence.

\begin{proposition}
    Let $A$ be an $R$-algebra. Then, the homotopy sheaves of
    $\Shaut_{\Mod_A}$ are
    \begin{equation*}
        \pi_i^s\Shaut_{\Mod_A}\iso\begin{cases}
            \ShAut_{\Mod_A}    &   \text{if $i=0$,}\\
            \THHB^0(A)^{\times}                 &   \text{if $i=1$,}\\
            \THHB^{1-i}(A)   &   \text{if $i\geq 2$,}
        \end{cases}
    \end{equation*}
    where $\ShAut_{\Mod_A}$ is the sheaf of groups with sections over $S$ the group
    $\Aut_{\Mod_{A\otimes_R S}}$, and $\THHB^*(A)$ is the Hochschild cohomology sheaf of $A$,
    which sends $S$ to $\THH^*_S(A\otimes_R S)$.
    \begin{proof}
        The description of $\pi_0^s\Shaut_{\Mod_A}$ is by definition. Since this is a sheaf
        of group-like $\EE_1$-spaces, the higher homotopy sheaves are independent of the
        basepoint chosen. The canonical basepoint is the identity functor $\id$, and it suffices to compute the
        homotopy sheaves of the loopsheaf $\Omega_{\id}\Shaut_{\Mod_A}$. Thus, one wants
        to compute the space of automorphisms of $\id:\Mod_A\rightarrow\Mod_A$ as a
        functor. This is nothing other than the space of automorphisms of $A$ as an
        $A^{\op}\otimes_R A$-module, which is precisely
        the space of units in the Hochschild cohomology algebra of $A$. The Hochschild
        cohomology algebra $\THHB^*(A)$ is a sheaf over $\Spec R$ because the category
        $\Mod_{A^{\op}\otimes_R A}$ satisfies \'etale hyperdescent.
    \end{proof}
\end{proposition}

This is a sheafy version of~\cite{toen-morita}*{Corollary 1.6}. When $X$ is quasi-compact and
quasi-separated, $\Mod_X\we\Mod_A$ for some $A$, so that the proposition also applies to the
automorphism sheaf of $\Mod_X$.

When $X$ is smooth, proper, and geometrically connected over $R$, $\THHB^0(X)^\times\iso\Gm$, since
$\THH^0_R(X)\iso\Hoh^0(X,\Oscr_X)$ by the Hodge spectral sequence for Hochschild
cohomology~\cite{swan}.  Moreover, if $R$ is an ordinary ring, and if $A$ is an ordinary
$R$-algebra or $X$ is an ordinary scheme, then the negative Hochschild cohomology groups
vanish, since projective resolutions exist.



The next theorem gives the main computational tool for determining $\Br^A(R)$. Throughout,
when writing $\ShBr^A$, it is assumed that $\Mod_A$ is chosen as the global basepoint of the
sheaf.

\begin{theorem}
    There is a fringed spectral sequence
    \begin{align*}
        \Eoh_2^{p,q}&=\begin{cases}
            \Hoh^p_{\et}(R,\pi_q^s\ShBr^A)  &   \text{if $q-p\geq 0$,}\\
            0                               &   \text{otherwise}
        \end{cases}\\
            &\Rightarrow\pi_{q-p}\ShBr^A(R),
    \end{align*}
    where
    \begin{equation*}
        \pi_i^s\ShBr^A\iso\begin{cases}
            0               &   \text{if $i=0$,}\\
            \ShAut_{\Mod_A}    &   \text{if $i=1$,}\\
            \THHB^0(A)^{\times}                 &   \text{if $i=2$,}\\
            \THHB^{2-i}(A)   &   \text{if $i\geq 3$,}
        \end{cases}
    \end{equation*}
    which converges completely if $A$ is smooth and proper or if $R$ and $A$ are ordinary rings.
    \begin{proof}
        The spectral sequence is nothing more than the descent spectral sequence of the
        previous section. The first statement about convergence follows because if $A$ is smooth
        and proper, the Hochschild cohomology of $A$ vanishes in sufficiently high degrees,
        so that the spectral sequence collapses after some finite stage. The second
        statement follows because, if $R$ and $A$ are ordinary, $\THH^{2-i}(A)=0$ for $i\geq 3$.
    \end{proof}
\end{theorem}

The theorem is especially strong when $R$ and $A$ are ordinary rings, or when $R$ is an
ordinary ring and one considers $\ShBr^X$ for a smooth, proper, geometrically connected $R$-scheme $X$. In either of
these cases the homotopy sheaves of the twisted Brauer sheaf are concentrated in two degrees,
$1$ and $2$. For instance,
\begin{equation*}
    \pi_i^s\ShBr^X\iso\begin{cases}
        0               &   \text{if $i=0$,}\\
        \ShAut_{\Mod_X}    &   \text{if $i=1$,}\\
        \Gm                 &   \text{if $i=2$,}\\
        0   &   \text{if $i\geq 3$.}
    \end{cases}
\end{equation*}
This means that the sheaf $\ShBr^X$ is an extension of Eilenberg-MacLane sheaves
\begin{equation*}
    K(\Gm,2)\rightarrow\ShBr^X\rightarrow K(\ShAut_{\Mod_X},1).
\end{equation*}
Since $\Gm$ is a sheaf of abelian groups, $K(\Gm,2)$ is an infinite loop space in
$\Shv_R^{\et}$. This implies that the sequence above can be delooped, and $\ShBr^X$ can be
identified as the fiber in the sequence
\begin{equation*}
    \ShBr^X\rightarrow K(\ShAut_{\Mod_X},1)\rightarrow K(\Gm,3).
\end{equation*}
Then, taking global sections, there is a fiber sequence
\begin{equation*}
    \ShBr^X(R)\rightarrow\map(\Spec R,K(\ShAut_{\Mod_X},1))\rightarrow\map(\Spec
    R,K(\Gm,3)).
\end{equation*}
We can point the spaces in this sequence by choosing the point $\Mod_X$ of
$\pi_0\ShBr^X(R)$. The spectral sequence degenerates into the long exact sequence of
homotopy groups associated to this fibration. In particular, there is an isomorphism
$\pi_2\ShBr^X(R)\iso\Gm(R)$ and an exact sequence
\begin{align*}
    0&\rightarrow\Hoh^1_{\et}(\Spec R,\Gm)\rightarrow\pi_1\ShBr^X(R)\rightarrow\ShAut_{\Mod_X}(R)\\
    &\rightarrow\Hoh^2_{\et}(\Spec
    R,\Gm)\rightarrow\pi_0\ShBr^X(R)\rightarrow\Hoh^1_{\et}(\Spec
    R,\ShAut_{\Mod_X})\rightarrow\Hoh^3_{\et}(\Spec R,\Gm).
\end{align*}
The meaning of exact here is just as in the beginning of the previous section. In
particular, there is an action of $\Hoh^2_{\et}(\Spec R,\Gm)$ on $\pi_0\ShBr^X(R)$, and the
fibers of $\pi_0\ShBr^X(R)\rightarrow\Hoh^1_{\et}(\Spec R,\ShAut_{\Mod_X})$ are precisely
the orbits of this action. The quotient $\Hoh^2_{\et}(\Spec R,\Gm)/\ShAut_{\Mod_X}(R)$ is
in bijection with the orbit of $\Mod_X$ in $\pi_0\ShBr^X(R)$.

The kernel of $\ShAut_{\Mod_X}(R)\rightarrow\Hoh^2_{\et}(\Spec R,\Gm)$ consists of those
elements that come from actual autoequivalences of $\Mod_X$.

An element of $\Hoh^1_{\et}(\Spec R,\ShAut_{\Mod_X})$ maps to $0$
in $\Hoh^3_{\et}(\Spec R,\Gm)$ if and only if it can be lifted to $\pi_0\ShBr^X(R)$. The
class in $\Hoh^3_{\et}(\Spec R,\Gm)$ represents the obstruction to lifting a cohomology
class in $\Hoh^1_{\et}(\Spec R,\ShAut_{\Mod_X})$ to an actual collection of gluing data to
obtain a twisted form of the stack $\StMod_X$. We will see that these obstructions
frequently vanish. This occurs when the gluing data can be made to act on an
object with less homotopical information, such as a scheme, as
opposed to the stable $\infty$-categories appearing in $\StMod_X$.

\subsection{The example of the introduction}\label{sec:curvesoverr}

Recall that $\Mod_{\PP^1_\RR}\we\Mod_{\RR Q}$ where $Q$ is quiver
$\bullet\rightrightarrows\bullet$. Since $\PP^1$ is Fano, the computation of Bondal and
Orlov~\cite{bondal-orlov} shows that $\ShAut_{\Mod_{\PP^1}}\iso\ZZ\times(\ZZ\rtimes\ShPGL_2)$.
Thus, by the vanishing of negative Hochschild cohomology for ordinary schemes,
the homotopy sheaves of $\ShBr^{\PP^1}$ are
\begin{equation*}
    \pi_i^s\ShBr^{\PP^1}=\begin{cases}
        0   &   \text{if $i=0$,}\\
        \ZZ\times(\ZZ\times\ShPGL_2)    &   \text{if $i=1$,}\\
        \Gm &   \text{if $i=2$,}\\
        0   &   \text{otherwise,}
    \end{cases}
\end{equation*}
where the degree $1$ term splits because $\ShPGL_2$ acts trivially on $\Pic(\PP^1)=\ZZ$.

In the descent spectral sequence for $\ShBr^{\PP^1}$ there is only one possible
non-zero differential, which is
$d_2:\ZZ\times(\ZZ\times\ShPGL_2(\RR))\rightarrow\Hoh^2_{\et}(\Spec\RR,\Gm)$.
But, it is clear that this is zero, because $\ZZ\times(\ZZ\times\ShPGL_2(\RR))$ survives to
the $\Eoh_\infty$-page to act as automorphisms of $\Mod_{\PP^1}$. Since $\Hoh^1_{\et}(\Spec k,\ZZ)=0$ for any field $k$, there is an
exact sequence of pointed sets
\begin{equation*}
    0\rightarrow\Br(\RR)\rightarrow\Br^{\PP^1}(\RR)\rightarrow\Hoh^1_{\et}(\Spec\RR,\ShPGL_2)\rightarrow *.
\end{equation*}
This sequence is not split, because the action of $\Br(\RR)$ on the
non-trivial point of $\Hoh^1_{\et}(\Spec \RR,\ShPGL_2)$ is trivial. However, since
$\Hoh^1_{\et}(\Spec\RR,\ShPGL_2)$ is the set of isomorphism classes of smooth projective genus
$0$ curves over $\RR$, the map $\Br^{\PP^1}(\RR)\rightarrow\Hoh^1_{\et}(\Spec\RR,\ShPGL_2)$
is indeed surjective, and the
set $\Br^{\PP^1}(\RR)$ consists of categories of twisted sheaves on genus $0$ curves.

We can compute the higher homotopy of $\ShBr^{\PP^1}(\RR)$ at the point $\Mod_{\PP^1}$. From
the spectral sequence,
\begin{equation*}
    \pi_i\ShBr^{\PP^1}(\RR)\iso\begin{cases}
        \ZZ\times(\ZZ\times\ShPGL_2(\RR))   &   \text{if $i=1$,}\\
        \RR^* &   \text{if $i=2$,}\\
        0   &   \text{if $i\geq 3$.}
    \end{cases}
\end{equation*}
Note that the fundamental group (at the point $\Mod_{\PP^1_\RR}$ is the automorphism group of $\Mod_{\PP^1_\RR}$.
The first $\ZZ$ is just translation, while the second corresponds to tensoring with $\Oscr(1)$.
The group $\pi_2$ is the group of invertible natural transformations between automorphisms.

The reader might be disturbed by an apparent asymmetry in the computation above. Namely,
what would happen if we did the calculation instead at the point $\Mod_C$ where $C$ is again
the curve $x^2+y^2+z^2=0$ in $\PP^2$ over $\RR$? In this case,
\begin{equation*}
     \pi_i^s\ShBr^C=\begin{cases}
        0   &   \text{if $i=0$,}\\
        \ZZ\times(\Roh p^1_*\mathds{G}_{m,C}\times\ShAut_C)    &   \text{if $i=1$,}\\
        \Gm &   \text{if $i=2$,}\\
        0   &   \text{otherwise.}
    \end{cases}
\end{equation*}
Here, $\ShAut_C$ is a form of $\ShPGL_2$ over $\RR$, and $p:C\rightarrow\Spec\RR$ is the structure map,
so $\Roh p^1_*\mathds{G}_{m,C}$ is the relative Picard sheaf. Then, by considering the Leray
spectral sequence for the sheaf $\mathds{G}_{m,C}$ and the map $p$, it is easy to see that
$\Pic(C)\rightarrow\Gamma(\Spec\RR,\Roh p^1_*\mathds{G}_{m,C})$ has cokernel equal to $\ZZ/2$.
Since we know that the Brauer group acts trivially, it follows that there is a non-zero
differential in the descent spectral sequence, and we obtain a filtration of pointed sets
\begin{equation*}
    0\rightarrow\Br^C(\RR)\rightarrow\Hoh^1_{\et}(\Spec\RR,\ShAut_C)\rightarrow *.
\end{equation*}
The difference between this computation and that for $\PP^1$ is simply because of the
dependence of the fiber on the basepoint for fibrations $X\rightarrow Y$, as discussed
in Section~\ref{sec:fringed}.

Thus, $\Mod_{\PP^1}$, $\Mod_{\PP^1}^{\HH}$, and $\Mod_C$ are the only
$3$ elements of $\Br^{\PP^1}$, which gives a positive answer to Question~\ref{q:grothendieck}.

\begin{theorem}
    Suppose that $A$ is an $\RR$-algebra such that $\CC\otimes_\RR A$ is derived Morita equivalent
    to $\PP^1$. Then, $A$ is derived Morita equivalent over $\RR$ to an ordinary $\RR$-algebra,
    either $\RR Q$, $\HH Q$, or the modulated quiver algebra associated to $C$.
\end{theorem}

\subsection{The stabilizer conjecture in the canonical (anti-)ample case}

The next theorem verifies the stabilizer conjecture in many cases.

\begin{theorem}\label{thm:stabilizerample}
    Let $X$ be a smooth, projective, and geometrically connected variety over a field $k$.
    Suppose that the canonical line bundle $\omega_X$ is either ample or anti-ample. Then,
    the stabilizer conjecture holds for $\Mod_X$. That is, the kernel of
    $\Br(k)\rightarrow\Br(X)$ is the same as the fiber of $\Br(k)\rightarrow\Br^X(k)$.
    \begin{proof}
        Consider the split exact sequence of sheaves of groups
        \begin{equation*}
            0\rightarrow\ZZ\times\ShPic_{X/k}\rightarrow\ShAut_{\Mod_X}\rightarrow\ShAut_X\rightarrow 0
        \end{equation*}
        given by the theorem of Bondal and Orlov~\cite{bondal-orlov}. The end of
        Section~\ref{subsec:ss} provides an exact sequence
        \begin{equation*}
            \pi_1\ShBr^X(k)\rightarrow\ShAut_{\Mod_X}(k)\rightarrow\Hoh^2_{\et}(k,\Gm)\rightarrow\pi_0\ShBr^X(k).
        \end{equation*}
        Thus, it suffices to show that the image of
        $\ShAut_{\Mod_X}(k)\rightarrow\Hoh^2_{\et}(k,\Gm)$ is precisely
        $\ker(\Br(k)\rightarrow\Br(X))$. By examining the exact sequence of sheaves above,
        it is clear that the only sections of $\ShAut_{\Mod_X}$ over $\Spec k$ that might
        not lift to automorphisms of $\Mod_X$ come from elements of $\ShPic_{X/k}(k)$ that do not
        lift to $\Pic(X)$. But, the cokernel of $\Pic(X)\rightarrow\ShPic_{X/k}(k)$
        injects into $\Br(k)$ as the kernel of $\Br(k)\rightarrow\Br(X)$ by the Leray
        spectral sequence. This completes the proof.
    \end{proof}
\end{theorem}

\section{Lifting Morita equivalences}\label{sec:moritadescent}

Let $k$ be a field, and let $A$ be a $k$-algebra. Up to this point, only algebras $B$ that
become derived Morita equivalent to $A$ after a finite separable extension $l/k$ have been
considered.

\begin{question}
    When is it the case that if $A$ and $B$ are $k$-algebras that are derived Morita equivalent over
    $\overline{k}$, then they are derived Morita equivalent over a finite separable extension of $k$.
\end{question}

This is a question about the smoothness of the stack of derived Morita equivalences between $A$ and
$B$. It is possible to solve it using the techniques of~\cite{ag} when $A$ is a smooth finite-dimensional hereditary $k$-algebra.
Recall that $A$ is hereditary if it has global dimension
$1$, and that $A$ is smooth if it has finite projective dimension over $A^{\op}\otimes_ RA$.

\begin{theorem}
    Let $A$ be a smooth finite-dimensional hereditary $k$-algebra. Then, if $B$ is derived Morita equivalent to $A$
    over $\overline{k}$, it is derived Morita equivalent to $A$ over some finite separable extension $l/k$.
    \begin{proof}
        Let $\ShMor_{A\rightarrow B}$ the sheaf of derived Morita equivalences from $A$ to $B$.
        Then, by hypothesis,
        $\ShMor_{A\rightarrow B}\rightarrow\Spec k$ is surjective on geometric points. Let $M$
        be a Morita equivalence over a field $l$. I can assume given $M$ that it is in fact
        a self-equivalence $\Mod_A\rightarrow\Mod_A$, and even the identity, viewed as a
        perfect complex of $\Mod_{A^{\op}\otimes_k A}$. But then, by~\cite{ag}*{Corollary
        5.9}, the cotangent complex of $\ShMor_{A\rightarrow B}$ at $M$ is equivalent to
        \begin{equation*}
            \Sigma^{-1}\End_{A^{\op}\otimes_k A}(A)^\vee.
        \end{equation*}
        The conditions on $A$ ensure
        that $\End_{A^{\op}\otimes_k A}(A)$ has homology (and
        hence, in this case, Tor-amplitude) contained in degrees $[-1,0]$.
        Since the base is a field, the dual $\End_{A^{\op}\otimes_k A}(A)^\vee$ has
        Tor-amplitude contained in $[0,1]$.  Thus, the cotangent complex has
        Tor-amplitude contained in degrees $[-1,0]$. Therefore, the sheaf of derived Morita
        equivalences is smooth when $A$ is a smooth finite-dimensional hereditary hereditary
        $k$-algebra. By
        Theorem~\cite{ag}*{Theorem 4.47}, it follows that there are \'etale local sections
        of $\ShMor_{A\rightarrow B}\rightarrow\Spec k$, as desired.
    \end{proof}
\end{theorem}

The theorem applies in particular to all path algebras, and hence also to the path algebras
of the quiver $\Omega_n$ from the previous section. There is also a global version of the
theorem, which has the same proof.

\begin{scholium}
    Suppose that $X$ is a regular noetherian scheme and that $\Ascr$ is a perfect sheaf of coherent
    algebras on $X$ such that $\Ascr_{k(x)}$ is smooth and hereditary for each point $x$ of $X$.
    If $\Bscr$ is another perfect sheaf of coherent algebras on $X$, and if
    $$\Mod_{\Bscr\otimes_{\Oscr_X}k(\overline{x})}\we\Mod_{\Ascr\otimes_{\Oscr_X} k(\overline{x})}$$ for each geometric
    point $\overline{x}$ of $X$, then there is an \'etale cover $U\rightarrow X$ such that
    $\Mod_{\Bscr\otimes_{\Oscr_X}\Oscr_U}\we\Mod_{\Ascr\otimes_{\Oscr_X}\Oscr_U}$.
\end{scholium}

\section{Examples}\label{sec:geoexamples}

The purpose of this section is to give a taste of the computational power of
the spectral sequence rather than to give a complete treatment in general. However, complete
computations are obtained for genus $0$ curves, quadric hypersurfaces, and twists of the
quiver $\Omega_n$. For curves of higher genus, only the outline of the theory is exposed, a
more detailed treatment being left to future work.

The spectral sequence makes it possible to describe representatives for all of the elements of $\pi_0\ShBr^A(k)$ or
$\pi_0\ShBr^X(k)$ in many cases. However, it is a much more difficult question to
decide when two representatives determine the same point in the set of connected components.
There are two reasons for this difficulty. First, it is in general a subtle problem to
determine the stabilizer of a point under the action of the Brauer group. In some good
cases, such as genus $0$ curves or quadrics, this is possible. But, for curves of higher
genus, for example, it is much harder to determine the stabilizer. The second problem is
that many sequences involve short exact sequences in nonabelian cohomology, where exactness
is only certain over basepoints.

Recall that I will abuse notation and write $\Mod_A$ and $\StMod^A$ interchangeably.

\subsection{Genus $0$ curves}

The reader can easily use the arguments in the introduction and
Section~\ref{sec:curvesoverr} to compute $\Br^{\PP^1}(k)$ for any field $k$.
There is a sequence
\begin{equation*}
    0\rightarrow\Br(k)\rightarrow\Br^{\PP^1}(k)\rightarrow\Hoh^1_{\et}(\Spec
    k,\ShPGL_2)\rightarrow *,
\end{equation*}
which is exact in the following senses. There is an action of $\Br(k)$ on $\Br^{\PP^1}(k)$,
which is faithful. Moreover, the action of $\Br(k)$ on the point $\Mod_{\PP^1}$ is free. The
map $\Br^{\PP^1}(k)\rightarrow\Hoh^1_{\et}(\Spec k,\ShPGL_2)$ is surjective, and two elements
of $\Br^{\PP^1}(k)$ lay over the same genus $0$ curve $C$ in $\Hoh^1_{\et}(\Spec k,\ShPGL_2)$
if and only they are in the same $\Br(k)$-orbit. A remark is in order about surjectivity,
as the end of Section~\ref{subsec:ss} implies that there is in general an obstruction.
However, it vanishes in this case as $\Hoh^1_{\et}(\Spec k,\ShPGL_2)$ is the
set of isomorphism of smooth projective genus $0$ curves over $k$. A class of
$\Br^{\PP^1}(k)$ mapping to $C\in\Hoh^1_{\et}(\Spec k,\ShPGL_2)$ can be constructed
explicitly be taking $\Mod_C$.

The interesting question is to determine
the stabilizers of the action of $\Br(k)$ on $\Mod_C$ for a genus $0$ curve without any
$k$-points. The curve $C$
is the Severi-Brauer variety of a unique degree $2$ central division algebra $D$ over $k$.
By Amitsur's theorem~\cite{gille-szamuely}*{Theorem 5.4.1}, the kernel of $\Br(k)\rightarrow\Br(k(C))$ is exactly the cyclic subgroup
generated by $[D]$. Since $\Br(C)\rightarrow\Br(k(C))$ is injective,
Theorem~\ref{thm:stabilizerample} says that the stabilizer is precisely
$([D])\subseteq\Br(k)$. It
follows that the orbit of $\Mod_C$ in $\Br^{\PP^1}(k)$ is in bijection with $\Br(k)/([D])$.

In summary, the noncommutative \'etale twists of $\PP^1$ are all determined by a genus $0$
curve $C$ and a Brauer class $\alpha\in\Br(k)$. The $\infty$-category of modules over this
noncommutative twist is $\Mod_C^{\alpha}$, the $\infty$-category of $\alpha$-twisted sheaves
on $C$. By using modulated quivers (for which, see~\cite{dlab-ringel}), all twists are derived Morita equivalent to ordinary
$k$-algebras, which answers Question~\ref{q:grothendieck} for the path algebra of
$\bullet\rightrightarrows\bullet$.

\subsection{Genus $1$ curves and modular representations}

Let $E$ be an elliptic curve over $k$. A group isomorphism $E\times E\rightarrow E\times E$
can be given by a matrix
\begin{equation*}
    f=\begin{pmatrix}
        f_1 &   f_2\\
        f_3 &   f_4
    \end{pmatrix}
\end{equation*}
of group isomorphisms $f_i:E\rightarrow E$. By using an isomorphism $E\iso{\hat{E}}$, where
$\hat{E}$ is the dual of $E$, one obtains
\begin{equation*}
    \tilde{f}=\begin{pmatrix}
        \hat{f}_4 &   -\hat{f}_2\\
        -\hat{f}_3 &   \hat{f}_1
    \end{pmatrix}.
\end{equation*}
Let $U(E)$ be the subgroup of group automorphisms $f:E\times E\rightarrow E\times E$
such that $\tilde{f}=f^{-1}$. Then, by Orlov~\cite{orlov-abelian}, there is an exact sequence
\begin{equation*}
    0\rightarrow\ZZ\times E\times\hat{E}\rightarrow\ShAut_{\Mod_E}\rightarrow U(E)\rightarrow 0.
\end{equation*}
From this, one can describe the elements of $\Br^E(k)$.

I will consider a special case, when
$U(E)\iso\SL_2(\ZZ)$, which happens for a non-CM elliptic curve. In this case, the
sequence reduces to
\begin{equation}\label{eq:2342}
    0\rightarrow\ZZ\times E\times\hat{E}\rightarrow\ShAut_{\Mod_E}\rightarrow\SL_2(\ZZ)\rightarrow 0.
\end{equation}
Let $\widetilde{\SL}_2(\ZZ)$ be the group generated by $x$, $y$, and $t$ with relations
$(xy)^3=t$, $y^4=t^2$, $xt=tx$, and $yt=ty$. Then, the quotient of $\widetilde{\SL}_2(\ZZ)$
by the central subgroup $(t)$ is isomorphic to $\SL_2(\ZZ)$. Moreover, there is a
homomorphism $\widetilde{\SL}_2(\ZZ)\rightarrow\ShAut_{\Mod_E}$ whose composition with
$\ShAut_{\Mod_E}\rightarrow\SL_2(\ZZ)$ is the surjection above. The element $t$ maps to the
translation functor. See~\cite{huybrechts}*{Section 9.5}. Since
$\ZZ\iso(t)\subseteq\widetilde{\SL}(\ZZ)$ is a central subgroup, it follows
from~\cite{serre-cohomologie}*{Proposition 42} that
\begin{equation*}
    \Hoh^1_{\et}(\Spec k,\widetilde{\SL}_2(\ZZ))\rightarrow\Hoh^1_{\et}(\Spec k,\SL_2(\ZZ))
\end{equation*}
is a bijection of pointed sets. Combining this fact with the exact sequence~\eqref{eq:2342}, one
easily proves the following lemma.

\begin{lemma}
    The natural map $\Hoh^1_{\et}(\Spec k,\ShAut_{\Mod_E})\rightarrow\Hoh^1_{\et}(\Spec
    k,\SL_2(\ZZ))$ is surjective.
\end{lemma}

Now, the descent spectral sequence for $\ShBr^E(k)$ yields an exact sequence
\begin{equation*}
    0\rightarrow\Br(k)\rightarrow\Br^E(k)\rightarrow\Hoh^1_{\et}(\Spec
    k,\ShAut_{\Mod_E})\rightarrow \Hoh^3_{\et}(\Spec k,\Gm).
\end{equation*}
(The exactness on the left follows from the stabilizer conjecture for $E$, which can be
proved by adapting the proof of the canonical (anti-)ample case to the non-CM elliptic curve
$E$ by using the explicit description of the sheaf of derived autoequivalences of $\Mod_E$.)
Let $v:\Br^E(k)\rightarrow\Hoh^1_{\et}(\Spec k,\ShAut_{\Mod_E})\rightarrow\Hoh^1_{\et}(\Spec k,\SL_2(\ZZ))$.
Since $\SL_2(\ZZ)$ is the constant \'etale sheaf, there is an equivalence
\begin{equation*}
    \Hoh^1_{\et}(\Spec k,\SL_2(\ZZ))\iso\Hom_{\cont}(\Gal(k),\SL_2(\ZZ)),
\end{equation*}
where $\Hom_{\cont}$ denotes continuous group homomorphisms.

\begin{proposition}
    To every twisted form $\Mrm$ of $\Mod_E$ there is a canonical modular representation of
    $\Gal(k)$.
\end{proposition}

Now, suppose that the $v$-invariant of $\Mrm$ is trivial. Then, using the exact sequence
\begin{equation*}
    \Hoh^1_{\et}(\Spec k,E\times\hat{E})\rightarrow\Hoh^1_{\et}(\Spec
    k,\ShAut_{\Mod_E})\rightarrow\Hoh^1_{\et}(\Spec k,\SL_2(\ZZ)),
\end{equation*}
it follows that the $\Br(k)$-orbit of $\Mrm$ corresponds to a class of $\Hoh^1_{\et}(\Spec
k,E\times E)$. The two copies of $E$ are not equal. One is $E$ acting on itself via
translations, the other is $\hat{E}$ acting on $\Mod_E$ by tensoring with line bundles.
The set $\Hoh^1_{\et}(\Spec k,E)$ contributes the categories $\Mod_C$ for homogeneous spaces
of $E$, while $\Hoh^1_{\et}(\Spec k,\hat{E})$ contributes the categories
$\Mod_E^\alpha$, where $\alpha\in\Br(E)/\Br(k)\subseteq\Hoh^1_{\et}(\Spec k,\hat{E})$, which
fits into the exact sequence
$0\rightarrow\Br(k)\rightarrow\Br(E)\rightarrow\Hoh^1_{\et}(k,\hat{E})\rightarrow\Hoh^3_{\et}(k,\Gm)$
coming from the Leray spectral sequence.

Consider twists $\Mrm\we\Mod_C$ where $C$ is a homogeneous space for $E$. It is impossible
at the moment to give a full treatment of the stabilizer of $\Mod_C$ in $\Br(k)$. The
same arguments used to prove the stabilizer conjecture in the canonical (anti-)ample case
can be used for a non-CM elliptic curve as well, which shows that the stabilizer is exactly
the kernel of $\Br(k)\rightarrow\Br(k(C))$. Until recently, very little was known about this kernel when
$C$ is a curve of genus higher than $0$. This has changed with the work of~\cite{han,
haile-han-wadsworth, ciperiani-krashen}. In~\cite{ciperiani-krashen}, the authors study this problem, and
show that for homogeneous spaces of curves over numbers fields or local fields, the kernel
can be computed algorithmically. They describe, for instance, a homogeneous space $C$ for
\begin{equation*}
    y^2+xy+y=x^3=x^2-10x-10
\end{equation*}
over $\QQ$
where the kernel, and hence stabilizer group, is isomorphic to $\ZZ/4\times\ZZ/2$. They also
show that for a homogeneous space over a local field or number field, the stabilizer is
finite~\cite{ciperiani-krashen}*{Proposition 4.11}. Over larger fields, they give an example
to show that the stabilizer need not be finite in general.
In the case of genus $0$ curves, the stabilizer is not only finite, but in all cases has
order at most $2$.


\subsection{Genus $g\geq 2$ curves}

Let $C$ be a smooth projective curve over $k$ having genus
$g\geq 2$. Then, $\omega_C$ is ample, so that the automorphism group of $\Mod_C$ can be
computed by Bondal and Orlov:
\begin{equation*}
    \Aut_{\Mod_C}\iso\ZZ\times(\Pic(C)\rtimes\Aut(C)).
\end{equation*}
Since there are no $\ZZ$-torsors over $\Spec k$, there is an exact sequence
\begin{equation*}
    \Hoh^1_{\et}(\Spec k,\ShPic^0_{C/k})\rightarrow\Hoh^1_{\et}(\Spec
    k,\ShAut_{\Mod_C})\rightarrow\Hoh^1_{\et}(\Spec k,\ShAut_C)\rightarrow\ast,
\end{equation*}
where $\ShPic^0_{C/k}$ is the Jacobian variety of $C$, and
where the map is surjective on the right since the surjection
$\ShAut_{\Mod_C}\rightarrow\ShAut_C$ splits.
There is again a sequence
\begin{equation*}
    \Br(k)\rightarrow\Br^{C}(k)\rightarrow\Hoh^1_{\et}(\Spec
    k,\ShPic^0_{C/k}\rtimes\ShAut_C)\rightarrow\Hoh^3_{\et}(\Spec k,\Gm),
\end{equation*}
with the same exactness properties as the sequence above for $\PP^1$. The kernel on the left
is precisely the kernel of $\Br(k)\rightarrow\Br(C)$ by Theorem~\ref{thm:stabilizerample}.

\begin{proposition}
    The twists $\Mrm$ of $\Mod_C$ are the categories $\Mod_D^\alpha$ for $D$ a twisted form of
    $C$ and $\alpha\in\Br_{\sep}(D)$, where
    $\Br_{\sep}(D)=\ker\left(\Br(D)\rightarrow\Br(D_{k^\sep})\right)$.
    \begin{proof}
        If $k$ is separable, then $\Br_{\sep}(D)=\Br(D)$.
        Since we can change the basepoint of $\Br^C(k)$ to $\Mod_D$ for any twisted form $D$ of
        $C$, it is enough to treat the classes in the fiber of
        \begin{equation*}
            \Br^C(k)\rightarrow\Hoh^1_{\et}(\Spec
            k,\ShAut_{\Mod_C})\rightarrow\Hoh^1_{\et}(\Spec k,\ShAut_C).
        \end{equation*}
        Write $F$ for the set of these points. Then, $F$ can be described by the exact sequence
        \begin{equation*}
            \Br(k)\rightarrow F\rightarrow\Hoh^1_{\et}(\Spec
            k,\ShPic^0_{C/k})\rightarrow\Hoh^3_{\et}(\Spec k,\Gm),
        \end{equation*}
        where the map $\Hoh^1_{\et}(\Spec k,\ShPic^0_{C/k})\rightarrow\Hoh^3_{\et}(\Spec
        k,\Gm)$ is induced from the sequence above.
        The Leray spectral sequence and the fact that $R^2p_*\mathds{G}_{m,C}=0$, where
        $p:C\rightarrow\Spec k$, shows that there is also an exact sequence
        \begin{equation*}
            \Br(k)\rightarrow\Br_{\sep}(C)\rightarrow\Hoh^1_{\et}(\Spec
            k,\ShPic^0_{C/k})\rightarrow\Hoh^3_{\et}(\Spec k,\Gm).
        \end{equation*}
        The map $\Br_{\sep}(C)\rightarrow F$ given by sending $\alpha\in\Br_{\sep}(C)$ to
        $\Mod_C^\alpha$ induces a map of these sequences, from which it follows that
        $\Br_{\sep}(C)$ surjects onto $F$, which proves the proposition.
    \end{proof}
\end{proposition}

The stabilizer of $\Mod_C$ is again the kernel of
$\Br(k)\rightarrow\Br(C)=\Br(k(C))$. As far as I know, when the genus is $g\geq 2$, almost nothing
is known about the stabilizer groups, except for the fact that it vanishes if $C$
has a $k$-point. In that case, the map $\Hoh^1_{\et}(\Spec
k,\ShPic^0_{C/k})\rightarrow\Hoh^3_{\et}(\Spec k,\Gm)$ is identically zero, because $\Hoh^3_{\et}(\Spec
k,\Gm)\rightarrow\Hoh^3_{\et}(C,\Gm)$ is injective (a $k$-point defines a section).

\subsection{Quadric hypersurfaces}

Assume for simplicity that the characteristic of $k$ is not $2$.

Consider the quadratic form
$q=x_0^2+\cdots+x_{2n-1}^2-x_{2n}^2$ on $k^{2n+1}$. Let
$X=X(q)$ be the quadric hypersurface in $\PP^{2n}$ cut out by $q$. I want to study
$\Br^X(k)$. As $X$ is Fano and $\Pic(X)=\ZZ$, the theorem of Bondal and
Orlov~\cite{bondal-orlov} says that $\Aut_{\Mod_X}=\ZZ\times(\ZZ\times\Aut(X))$.
Therefore, every element of $\Br^X(k)$ is $\Mod_Y^\alpha$ for a twisted form $Y$ of $X$ and
a Brauer class $\alpha\in\Br(k)$. Every such $Y$ is determined by another nondegenerate
quadratic form $p$ on $k^{2n+1}$. The interesting question is then what is the stabilizer of
$\Mod_Y$ for such a twist $Y$; in other words, by Theorem~\ref{thm:stabilizerample}, what is the kernel of
$\Br(k)\rightarrow\Br(k(Y))$? This is in fact a classical question.
Let $C_0(p)$ denote the even Clifford algebra of $p$, which
is a central simple algebra. To compute the kernel of $\Br(k)\rightarrow\Br(Y)$, it is
enough to compute the kernel of $\Br(k)\rightarrow\Br(k(Y))$, since $Y$ is smooth. A
division algebra $D$ is in the kernel if and only if the index of $D\otimes_k k(Y)$ is $1$.
But, this index was computed to be
\begin{equation*}
    \gcd\{\ind(D),2^{n-1}\ind(D\otimes_k C_0(p))\},
\end{equation*}
by, for instance, Merkurjev-Panin-Wadsworth~\cite{merkurjev-panin-wadsworth}. Because
$C_0(p)$ has degree a power of $2$, the kernel must be $2$-primary. Therefore, if
$n>1$, the kernel is always $0$. The case $n=1$ was already handled in the case of genus $0$
curves. The following theorem summarizes the situation. Note that the statement about
surjectivity follows for the same reason as for genus $0$ curves; namely, concrete models
can be constructed by taking a twist $Y\in\Hoh^1_{\et}(\Spec k,\ShPSO(q))$ and then
considering $\Mod_Y\in\Br^X(k)$.

\begin{theorem}
    Suppose that $n>1$ and that $X$ is the quadric hypersurface in $\PP^{2n}$
    considered above. Then, there is a sequence
    \begin{equation*}
        0\rightarrow\Br(k)\rightarrow\Br^X(k)\rightarrow\Hoh^1_{\et}(\Spec
        k,\ShPSO(q))\rightarrow *,
    \end{equation*}
    which is exact in the sense that the action of $\Br(k)$ on $\Br^X(k)$ is free, and two
    elements of $\Br^X(k)$ map to the same element of $\Hoh^1_{\et}(\Spec k,\ShPSO(q))$ if and
    only if they are in the same $\Br(k)$-orbit.
\end{theorem}

Now, consider $q=x_0^2+\cdots+x_{2n-2}^2-x_{2n-1}^2$ on $k^{2n}$, and let
$X=X(q)$ be the quadric hypersurface in $\PP^{2n-1}$ cut out by $q$. Again, in this case every
class of $\Br^X(k)$ is $\Mod_Y^{\alpha}$ where $Y$ is a twist of $X$ (an involution variety) and $\alpha\in\Br(k)$.
To consider the stabilizer of $\Br(k)$ on $Y$, it suffices to compute the index of
$D\otimes_k k(Y)$ as above. By~\cite{merkurjev-panin-wadsworth}, this is
\begin{equation*}
    \ind(D\otimes_k k(Y))=\gcd\{\ind(D),2^{n-2}\ind(D\otimes_k C(p)\},
\end{equation*}
where $C(p)$ is the full Clifford algebra of $p$. The Clifford algebra $C(p)$ has index a
$2$-power, so that the kernel is once again $2$-primary. Therefore, if $n>2$, the kernel vanishes.
When $n=2$, the stabilizer of the quadric hypersurface $Y$ in $\Br^X(k)$ is generated by the
central simple algebra $C(p)$.

\begin{theorem}
    Suppose that $n\geq 1$, and let $X$ be the quadric hypersurface in
    $\PP^{2n-1}$ considered above. Then, there is a sequence
    \begin{equation*}
        0\rightarrow\Br(k)\rightarrow\Br^X(k)\rightarrow\Hoh^1_{\et}(\Spec k,\ShPSO(q))\rightarrow *,
    \end{equation*}
    which is exact in the sense that the action of $\Br(k)$ on $\Br^X(k)$ is faithful, and
    two elements of  $\Br^X(k)$ map to the same element of $\Hoh^1_{\et}(\Spec k,\ShPSO(q))$ if and
    only if they are in the same $\Br(k)$-orbit. Moreover, if $n>2$, the action is free. If
    $n=2$, then the stabilizer of the quadric surface $Y$ associated to a nondegenerate
    quadratic form $p$ is generated by the Clifford algebra $C(p)$.
\end{theorem}

\begin{remark}
    Using the work of Merkurjev, Panin, and Wadsworth, the same game can be played for the
    twisted flag varieties associated to any classical semisimple adjoint linear algebraic group.
\end{remark}

If $X$ is a quadric hypersurface, then $\Mod_X\we\Mod_A$ for an ordinary associative algebra
$A$, as follows from a theorem of Kapranov~\cite{kapranov}. The classification theorem says that every
object of $\Br^X(k)$ is equivalent to $\Mod_B$ for some ordinary $k$-algebra $B$, giving
another positive answer to Question~\ref{q:grothendieck}.

\subsection{Noncommutative Severi-Brauer varieties}\label{sec:ncsbv}

Kontsevich and Rosenberg~\cite{kontsevich-rosenberg} introduced a noncommutative space
$\NP^{n-1}$ that represents the functor which takes an associative algebra $A$ to the set of
quotients of $A^n$ that are locally isomorphic to $A$ in a flat topology on associative
rings. They thus called it the noncommutative projective space.
They identified its derived category with the
derived category of finite representations of the quiver $\Omega_n$.
\begin{figure}[h]
    \centering
    \begin{tikzpicture}
        \draw [fill] (-.1,.05) circle [radius=.05] node[left] {$a$};
        \draw [fill] (1.1,.05) circle [radius=.05] node[right] {$b$};
        \draw [->] (0,.1) to [out=45, in=135] node[above] {$s_1$} (1,.1);
        \draw [->] (0,.05) to [out=30, in=150] (1,.05);
        \draw [->] (0,0) to [out=-30, in=210] node[below] {$s_n$} (1,0) ;
        \draw [dotted] (.5,.17) to (.5,-.13);
    \end{tikzpicture}
    \caption{The quiver $\Omega_n$.}
\end{figure}
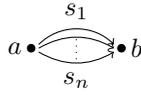
Thus, I consider
$\Mod_{k\Omega_n}$ as the model for noncommutative projective space. Except when $n=2$, this is
\emph{not} the derived category of $\PP^{n-1}$. Nevertheless,
in~\cite{miyachi-yekutieli}, Miyachi and Yekutieli showed another similarity between
$\NP^{n-1}$ and $\PP^{n-1}$ by computing the group of equivalences of
$\Mod_{k\Omega_n}$ and showing that it is $\ZZ\times(\ZZ\rtimes\ShPGL_n(k))$. It follows that
for every $\PGL_n$-torsor $P$ over $k$, there is a well-defined twisted form of
$\Mod_{k\Omega_n}$, which I will denote $\Mrm^P$. But, the $\PGL_n$-torsors are in
one-to-one correspondence with Severi-Brauer varieties. So, $\Mrm^P$ is a
noncommutative twist of the Severi-Brauer variety $P$.

\begin{theorem}
    There is a bijection between the $\Br(k)$-sets $\Br^{\PP^n}(k)$ and $\Br^{\NP^n}(k)$.
\end{theorem}

Once again, these can be described using path algebras for modulated quivers, so there is a
positive answer to Question~\ref{q:grothendieck}.

%

%


\begin{bibdiv}
\begin{biblist}
    


\bib{ag}{article}{
    author = {Antieau, Benjamin},
    author = {Gepner, David},
    title = {Brauer groups and \'etale cohomology in derived algebraic geometry},
    journal = {ArXiv e-prints},
    eprint = {http://arxiv.org/abs/1210.0290},
    year = {2012},
}











\bib{beilinson}{article}{
    author={Be{\u\i}linson, A. A.},
    title={Coherent sheaves on ${\bf P}^{n}$ and problems in linear
    algebra},
    journal={Funktsional. Anal. i Prilozhen.},
    volume={12},
    date={1978},
    number={3},
    pages={68--69},
    issn={0374-1990},
}

\bib{bzfn}{article}{
    author={Ben-Zvi, David},
    author={Francis, John},
    author={Nadler, David},
    title={Integral transforms and Drinfeld centers in derived algebraic
    geometry},
    journal={J. Amer. Math. Soc.},
    volume={23},
    date={2010},
    number={4},
    pages={909--966},
    issn={0894-0347},
}





\bib{bondal-vandenbergh}{article}{
    author={Bondal, A.},
    author={van den Bergh, M.},
    title={Generators and representability of functors in commutative and noncommutative geometry},
    journal={Mosc. Math. J.},
    volume={3},
    date={2003},
    number={1},
    pages={1--36, 258},
    issn={1609-3321},
}

\bib{bondal-orlov}{article}{
    author={Bondal, Alexei},
    author={Orlov, Dmitri},
    title={Reconstruction of a variety from the derived category and groups of autoequivalences},
    journal={Compositio Math.},
    volume={125},
    date={2001},
    number={3},
    pages={327--344},
    issn={0010-437X},
}

\bib{bondal-orlov-coherent}{article}{
    author={Bondal, A.},
    author={Orlov, D.},
    title={Derived categories of coherent sheaves},
    conference={
        title={Proceedings of the ICM},
        address={Beijing},
        date={2002},
    },
    book={
    publisher={Higher Ed. Press},
    place={Beijing},
    },
    date={2002},
    pages={47--56},
}
 
 
\bib{bousfield-kan}{book}{
    author={Bousfield, A. K.},
    author={Kan, D.  M.},
    title={Homotopy limits, completions and localizations},
    series={Lecture Notes in Mathematics, Vol. 304},
    publisher={Springer-Verlag},
    place={Berlin},
    date={1972},
    pages={v+348},
}
 
\bib{brown-gersten}{article}{
    author={Brown, Kenneth S.},
    author={Gersten, Stephen M.},
    title={Algebraic $K$-theory as generalized sheaf cohomology},
    conference={title={Algebraic K-theory, I: Higher K-theories (Proc. Conf., Battelle Memorial Inst., Seattle, Wash., 1972)}, },
    book={
        publisher={Springer},
        place={Berlin},
    },
    date={1973},
    pages={266--292. Lecture Notes in Math., Vol.  341},
}

\bib{caldararu}{article}{
    author = {C\u{a}ld\u{a}raru, Andrei},
    title = {Derived categories of twisted sheaves on Calabi-Yau manifolds},
    note = {Ph.D. thesis, Cornell University (2000)},
    eprint = {http://www.math.wisc.edu/~andreic/},
}

\bib{ciperiani-krashen}{article}{
    author={Ciperiani, Mirela},
    author={Krashen, Daniel},
    title={Relative Brauer groups of genus 1 curves},
    journal={Israel J. Math.},
    volume={192},
    date={2012},
    pages={921--949},
    issn={0021-2172},
}




\bib{cline-parshall-scott}{article}{
    author={Cline, E.},
    author={Parshall, B.},
    author={Scott, L.},
    title={Derived categories and Morita theory},
    journal={J. Algebra},
    volume={104},
    date={1986},
    number={2},
    pages={397--409},
    issn={0021-8693},
}

\bib{dejong}{article}{
    author={de Jong, Aise Johan},
    title={A result of Gabber},
    eprint={http://www.math.columbia.edu/~dejong/},
}





\bib{dlab-ringel}{article}{
    author={Dlab, Vlastimil},
    author={Ringel, Claus Michael},
    title={On Algebras of Finite Representation Type},
    journal={J. Algebra},
    volume={33},
    year={1975},
    pages={306-394},
}

%
 





 

\bib{gille-szamuely}{book}{
    author={Gille, Philippe},
    author={Szamuely, Tam{\'a}s},
    title={Central simple algebras and Galois cohomology},
    series={Cambridge Studies in Advanced Mathematics},
    volume={101},
    publisher={Cambridge University Press},
    place={Cambridge},
    date={2006},
    pages={xii+343},
}

\bib{ginzburg}{article}{
    author = {Ginzburg, Victor},
    title = {Lectures on noncommutative geometry},
    journal = {ArXiv e-prints},
    eprint = {http://arxiv.org/abs/math/0506603},
    year = {2005},
}
 


%
%


 
\bib{grothendieck-brauer-3}{incollection}{
    author={Grothendieck, Alexander},
    title={Le groupe de Brauer. III. Exemples et compl\'ements},
    conference={
        title={Dix Expos\'es sur la Cohomologie des Sch\'emas},
    },
    book={
        publisher={North-Holland},
        place={Amsterdam},
    },
    date={1968},
    pages={88--188},
}

\bib{haile-han-wadsworth}{article}{
    author={Haile, Darrell E.},
    author={Han, Ilseop},
    author={Wadsworth, Adrian R.},
    title={Curves $C$ that are cyclic twists of $Y^2=X^3+c$ and the relative Brauer groups
    $Br(k(C)/k)$},
    journal={Trans. Amer. Math. Soc.},
    volume={364},
    date={2012},
    number={9},
    pages={4875--4908},
    issn={0002-9947},
}

\bib{han}{article}{
    author={Han, Ilseop},
    title={Relative Brauer groups of function fields of curves of genus one},
    journal={Comm. Algebra},
    volume={31},
    date={2003},
    number={9},
    pages={4301--4328},
    issn={0092-7872},
}

\bib{happel}{article}{
    author={Happel, Dieter},
    title={On the derived category of a finite-dimensional algebra},
    journal={Comment. Math. Helv.},
    volume={62},
    date={1987},
    number={3},
    pages={339--389},
    issn={0010-2571},
}



\bib{huybrechts}{book}{
    author={Huybrechts, D.},
    title={Fourier-Mukai transforms in algebraic geometry},
    series={Oxford Mathematical Monographs},
    publisher={The Clarendon Press Oxford University Press},
    place={Oxford},
    date={2006},
    pages={viii+307},
    isbn={978-0-19-929686-6},
    isbn={0-19-929686-3},
}


\bib{kapranov}{article}{
    author={Kapranov, M. M.},
    title={On the derived categories of coherent sheaves on some homogeneous
    spaces},
    journal={Invent. Math.},
    volume={92},
    date={1988},
    number={3},
    pages={479--508},
    issn={0020-9910},
}



\bib{keller}{article}{
    author={Keller, Bernhard},
    title={Deriving DG categories},
    journal={Ann. Sci. \'Ecole Norm. Sup. (4)},
    volume={27},
    date={1994},
    number={1},
    pages={63--102},
    issn={0012-9593},
}

\bib{kontsevich-rosenberg}{article}{
    author={Kontsevich, Maxim},
    author={Rosenberg, Alexander L.},
    title={Noncommutative smooth spaces},
    conference={
    title={The Gelfand Mathematical Seminars, 1996--1999},
    },
    book={
    series={Gelfand Math. Sem.},
    publisher={Birkh\"auser Boston},
    place={Boston, MA},
    },
    date={2000},
    pages={85--108},
}



%

\bib{lieblich}{article}{
    author={Lieblich, Max},
    title={Twisted sheaves and the period-index problem},
    journal={Compositio Math.},
    volume={144},
    date={2008},
    number={1},
    pages={1--31},
    issn={0010-437X},
}

\bib{htt}{book}{
      author={Lurie, Jacob},
       title={Higher topos theory},
      series={Annals of Mathematics Studies},
   publisher={Princeton University Press},
     address={Princeton, NJ},
        date={2009},
      volume={170},
        ISBN={978-0-691-14049-0; 0-691-14049-9},
}



\bib{dag11}{article}{
    author={Lurie, Jacob},
    title={Derived algebraic geometry XI: descent theorems},
    date={2011},
    eprint={http://www.math.harvard.edu/~lurie/},
}


\bib{ha}{article}{
    author={Lurie, Jacob},
    title={Higher algebra},
    date={2012},
    eprint={http://www.math.harvard.edu/~lurie/},
}


\bib{merkurjev-panin-wadsworth}{article}{
    author={Merkurjev, A. S.},
    author={Panin, I. A.},
    author={Wadsworth, A. R.},
    title={Index reduction formulas for twisted flag varieties. I},
    journal={$K$-Theory},
    volume={10},
    date={1996},
    number={6},
    pages={517--596},
    issn={0920-3036},
}

\bib{miyachi-yekutieli}{article}{
    author={Miyachi, Jun-ichi},
    author={Yekutieli, Amnon},
    title={Derived Picard groups of finite-dimensional hereditary algebras},
    journal={Compositio Math.},
    volume={129},
    date={2001},
    number={3},
    pages={341--368},
    issn={0010-437X},
}



\bib{orlov-abelian}{article}{
    author={Orlov, D. O.},
    title={Derived categories of coherent sheaves on abelian varieties and equivalences between them},
    journal={Izv. Ross. Akad. Nauk Ser. Mat.},
    volume={66},
    date={2002},
    number={3},
    pages={131--158},
    issn={0373-2436},
    translation={
        journal={Izv. Math.},
        volume={66},
        date={2002},
        number={3},
        pages={569--594},
        issn={1064-5632},
    },
}

\bib{quillen}{article}{
    author={Quillen, Daniel},
    title={Higher algebraic $K$-theory. I},
    conference={
        title={Algebraic $K$-theory, I: Higher $K$-theories (Proc. Conf.,
        Battelle Memorial Inst., Seattle, Wash., 1972)},
    },
    book={
        publisher={Springer},
        place={Berlin},
    },
    date={1973},
    pages={85--147. Lecture Notes in Math., Vol.  341},
}

\bib{rickard}{article}{
    author={Rickard, Jeremy},
    title={Morita theory for derived categories},
    journal={J. London Math. Soc. (2)},
    volume={39},
    date={1989},
    number={3},
    pages={436--456},
    issn={0024-6107},
}

\bib{rickard-equivalences}{article}{
    author={Rickard, Jeremy},
    title={Derived equivalences as derived functors},
    journal={J. London Math. Soc. (2)},
    volume={43},
    date={1991},
    number={1},
    pages={37--48},
    issn={0024-6107},
}


\bib{rouquier-zimmermann}{article}{
    author={Rouquier, Rapha{\"e}l},
    author={Zimmermann, Alexander},
    title={Picard groups for derived module categories},
    journal={Proc. London Math. Soc. (3)},
    volume={87},
    date={2003},
    number={1},
    pages={197--225},
    issn={0024-6115},
}

%
%
 
\bib{schwede-shipley}{article}{
author={Schwede, Stefan},
author={Shipley, Brooke},
title={Stable model categories are categories of modules},
journal={Topology},
volume={42},
date={2003},
number={1},
pages={103--153},
issn={0040-9383},
}

\bib{serre-cohomologie}{book}{
    author={Serre, Jean-Pierre},
    title={Cohomologie galoisienne},
    series={Lecture Notes in Mathematics},
    volume={5},
    edition={5},
    publisher={Springer-Verlag},
    place={Berlin},
    date={1994},
    pages={x+181},
    isbn={3-540-58002-6},
}

 

\bib{swan}{article}{
    author={Swan, Richard G.},
    title={Hochschild cohomology of quasiprojective schemes},
    journal={J. Pure Appl. Algebra},
    volume={110},
    date={1996},
    number={1},
    pages={57--80},
    issn={0022-4049},
}

 
%
%

\bib{toen-morita}{article}{
    author={To{\"e}n, Bertrand},
    title={The homotopy theory of $dg$-categories and derived Morita theory},
    journal={Invent. Math.},
    volume={167},
    date={2007},
    number={3},
    pages={615--667},
    issn={0020-9910},
}

   

\bib{toen-derived}{article}{
    author = {To{\"e}n, Bertrand},
    title = {Derived Azumaya algebras and generators for twisted derived categories},
    journal = {Invent. Math.},
    issn = {0020-9910},
    volume = {189},
    number = {3},
    year = {2012},
    pages = {581-652},
    url = {http://dx.doi.org/10.1007/s00222-011-0372-1},
}


\bib{toen-vaquie}{article}{
    author={To{\"e}n, Bertrand},
    author={Vaqui{\'e}, Michel},
    title={Moduli of objects in dg-categories},
    journal={Ann. Sci. \'Ecole Norm. Sup. (4)},
    volume={40},
    date={2007},
    number={3},
    pages={387--444},
    issn={0012-9593},
}



\bib{vandenberg-blowing}{article}{
    author={Van den Bergh, Michel},
    title={Blowing up of non-commutative smooth surfaces},
    journal={Mem. Amer. Math. Soc.},
    volume={154},
    date={2001},
    number={734},
    pages={x+140},
    issn={0065-9266},
}



\bib{wang}{article}{
    author={Wang, Shianghaw},
    title={On the commutator group of a simple algebra},
    journal={Amer. J. Math.},
    volume={72},
    date={1950},
    pages={323--334},
    issn={0002-9327},
}

\end{biblist}
\end{bibdiv}
\end{document}